\documentclass{article}
\usepackage[a4paper,  top=1in, bottom=1.5in, left=1.25in, right=1.25in]{geometry}
\usepackage[utf8]{inputenc}
\usepackage{lmodern}
\usepackage{amsmath}
\usepackage{amssymb}
\usepackage{graphicx}
\usepackage[algoruled]{algorithm2e}
\usepackage{listings}
\usepackage{caption}
\usepackage{color}
\usepackage{authblk}
\usepackage{array} 
\usepackage{epstopdf}
\newcolumntype{L}{>{\displaystyle}l} 
\newcolumntype{C}{>{\displaystyle}c} 
\newcolumntype{R}{>{\displaystyle}r} 
\usepackage{parskip}

\usepackage{hyperref}

\def\fatLambda{\boldsymbol{\Lambda}}

\newcommand{\veps}{\varepsilon}
\newcommand{\beq}{\begin{equation}}
\newcommand{\eeq}{\end{equation}}

\def\fatA{\mathbf{A}}
\def\fatB{\mathbf B}

\def\fatc{\mathbf c}

\def\fate{\mathbf e}
\def\fatE{\mathbf E}
\def\fatf{\mathbf f}
\def\fatF{\mathbf F}

\def\fatH{\mathbf H}

\def\fatI{\mathbf I}
\def\fatJ{\mathbf J}
\def\fatm{\mathbf m}
\def\fatn{\mathbf n}
\def\fatO{\mathbf O}
\def\fatQ{\mathbf Q}
\def\fatr{\mathbf r}
\def\dfatr{\delta\mathbf r}
\def\hfatr{\hat{\fatr}}

\def\fatx{\mathbf x}
\def\tfatx{\tilde{\fatx}}
\def\dfatx{\delta\fatx}
\def\Dfatx{\Delta\fatx}
\def\faty{\mathbf y}
\def\fatu{\mathbf u}

\def\fatv{\mathbf v}

\def\fatzero{\mathbf 0}
\def\fat1{\mathbf 1}
\def\WbE{W_{\rm SRBE}}

\def\WbEtot{W_{\rm SRBE,tot}}

\def\WfE{W_{\rm SRFE}}
\def\WfEtot{W_{\rm SRFE,tot}}
\def\bbR{\mathbb R}
\def\Dt{\Delta t}
\def\DTs{\Delta \mathcal{T}}
\def\DT{\Delta T}
\def\Expc{\mathbb{E}}
\def\Tacc{T_{\rm acc}}
\def\Dtau{\Delta \tau} 
\def\Dtaus{\Delta \tau_{\rm s}}

\def\fatrho{\boldsymbol{\rho}}
\def\hfatrho{\hat{\fatrho}}

\def\dr{\delta r}

\def\dx{\delta x}

\def\tr{{\rm trans}}

\def\tfatA{\tilde{\fatA}}

\def\fateta{\boldsymbol{\eta}}
\def\fattheta{\boldsymbol{\theta}}

\def\fatPhi{\boldsymbol{\Phi}}
\def\ninf{\|_\infty}

\def\calB{\mathcal{B}}
\def\calC{\mathcal{C}}
\def\calD{\mathcal{D}}
\def\calF{\mathcal{F}}
\def\calI{\mathcal{I}}
\def\calJ{\mathcal{J}}
\def\calK{\mathcal{K}}

\def\calO{\mathcal{O}}

\def\d{\text{d}}

\begin{document}

\title{Adaptive time integration of mechanical forces in center-based models for
biological cell populations}
\author{Per Lötstedt, Sonja Mathias\\ Department of Information Technology. \\  Uppsala University, SE-751 05 Uppsala, Sweden
}
\date{\today}


\maketitle

\begin{abstract}
Center-based models are used to simulate the mechanical behavior of biological cells during embryonic development or cancer growth. 
To allow for the simulation of biological populations potentially growing from a few individual cells to many thousands or more, these models have to be numerically efficient, while being reasonably accurate on the level of individual cell trajectories. In this work, we increase the robustness, accuracy, and efficiency of the simulation of
center-based models
by choosing the time steps adaptively in the numerical method. 
We investigate the gain in using single rate time stepping for the forward and backward Euler methods, based on local estimates of the numerical errors and the stability of the method in the case of the explicit forward Euler method. 
Furthermore, we propose a multirate time stepping scheme that simulates regions with high local force gradients (e.g.\ as they happen after cell division) with multiple smaller time steps within a larger single time step for regions with smoother forces. These methods are compared for different model systems in numerical experiments. 
We conclude that the adaptive single rate forward Euler method results in significant gains in terms of reduced wall clock times for the simulation of a linearly growing tissue, while at the same time eliminating the need for manual determination of a suitable time step size.
\end{abstract}

{\bf Keywords}: {center-based model, time integration, adaptivity, gradient system,  numerical methods}

{\bf Mathematics Subject Classification (2020)}: 65Z05, 92C15,  92-10
\section{Introduction}

In computational models of tissue mechanics, biological cells grow, divide, migrate, and die. 
They are exposed to mechanical forces from other cells and the environment during their life time. In off-lattice cell-based models, each cell is represented individually and moves due to these forces in continuous two or three dimensional space. Off-lattice cell-based models are used in simulations to study e.g.\ tumours \cite{VLPJD15, Lima21, LFJCLMWC10}, morphogenesis \cite{GKV19}, and colonies of bacteria \cite{Hellweger16}.
The governing equations for the motion of the aggregation of cells constitute a system of ordinary differential equations (ODEs) derived from Newton's second law \cite{OptABM, ByDr09,  VLPJD15, MCBH20}. 
These ODEs are satisfied by the coordinates of the cell centers in a cell-centered or center-based model (CBM) and the forces between the cells depend on their relative distance \cite{DH05}. Another type of off-lattice model based on the time integration of such an ODE system is a vertex-based model (VBM) \cite{VLPJD15}. There the forces are applied in the vertices of the cell boundary. 
In either case the system needs to be solved numerically
at discrete time points. 

In addition to the cell movement cell-based models incorporate the aforementioned individual cell behaviors such as proliferation, migration and apoptosis which affect the mechanics of the population. As a result, there are different time scales in the ODE system. After cell proliferation or division for example, force magnitudes are large and short time steps are necessary in the numerical solution to resolve a fast time scale. Long time steps are possible between cell divisions when the cells move on a slow scale.
However, in a fixed time stepping scheme the time step size is constant and hence dictated by the minimal step size required after proliferation. This makes it very inefficient if the population behavior is highly dynamic. In fact, adaptive time stepping is mentioned as one of the challenges in the simulation of multicellular tissues in a recent survey \cite{FlOs20}. To make things worse, choosing too large time step sizes potentially results in numerically instable solutions when using explicit schemes. And even if the time step ensures stability, cell trajectories of daughter cells after division may be physically incorrect leading to geometrical differences on the population level indistinguishable from parameter value effects unless the time step size is reduced further \cite{MCBH20}. As such the correct choice of time step size is very important for the model behavior and consequently for the validity of model conclusions.
The effects of constant time steps on the quality of the solutions of a VBM are reported in \cite{KuBaFl17}.

A simple way of ensuring numerical stability in CBMs is adjusting the time step size based on a threshold on the spatial displacement of the cells between time steps \cite{Atwell2016, KBMH, VanLiedekerke2018, SchMH05}. 
This threshold is usually chosen significantly smaller than the cellular radius. If it is violated, the time step is decreased and the cell positions are recalculated. 
Additionally, one can also define a minimum spatial step size which if fulfilled will lead to the time step being increased \cite{VanLiedekerke2018}. 
The implementation of this simple heuristic approach (even without a minimal spatial displacement) can already lead to significant reduction in computational cost, as discussed in \cite{Atwell2016}. 
Two time steps are determined in \cite{KBMH} where a short step used to integrate the rapidly moving cells and the majority of the cells are advanced in time by a longer step. 
By limiting the spatial displacement in each time step, the drawback is that the modeler needs to choose the exact threshold value for this displacement. 
The correct values may need to be determined by trial-and-error until the simulation results no longer show any change. 
Another method for CBMs avoiding the need for the spatial displacement threshold was proposed in \cite{Atwell2016} based on an embedded Runge-Kutta scheme, in particular the Dormand-Prince 853 scheme \cite[p. 181]{HaNoWa}. 
Here an order eight Runge-Kutta method was used to calculate the solution itself, along with two lower order methods to provide an estimate of the numerical error. Such a high order, however, may not be necessary for cell simulations and its computational cost may be prohibitive.


In this paper, we compare methods to choose the time step in the numerical solution of the ODE system governing the motion of CBMs. The methods we propose vary the time step size not by checking a threshold, but by directly calculating the time step size necessary to satisfy a desired accuracy on the level of the cell trajectories. 
As such, the time step size is determined by the properties of the ODE system and the numerical method instead of the physical properties of the cells or geometric conditions. 
More specifically, the numerical solution is advanced by the backward and forward Euler methods with time steps such that the local estimate of the discretization error is bounded by a given parameter. 
In addition, there is a bound on the time step in the forward Euler method for the solution to remain stable. 
As a result the need to guess a constant time step for the whole time interval of interest is eliminated and at the same time the numerical errors are controlled by changing the time step sizes, thereby increasing the robustness of center-based model simulations.

The forward and backward Euler methods are of first order accuracy and the former is commonly used in CBMs \cite{Delile2017MecaGen, osborne2017comparing}, both due to its simplicity and the fact that for cell-based models the modeling errors usually dominate and hence high order accuracy of the ODE solver is not required. 
This is also the conclusion drawn in \cite{KuBaFl17}.
The backward Euler method is an implicit method and as such associated with a higher computational cost per time step than the explicit forward Euler method. This prohibits its use in combination with a fixed time stepping for many practical experimental setups \cite{Atwell2016, mathias2022cbmos}. However, its improved stability properties are suitable for stiff systems as they promise larger step sizes at mechanical equilibrium. We will  investigate whether the use of an adaptive time step size can render the backward Euler method computationally beneficial as well. 

Based on these methods, we consider two different types of adaptive schemes. 
In a globally adaptive single rate method, the same time step is used for all cells.
In a locally adaptive multirate method, on the other hand, the time step is shorter for cells subject to a rapid change and longer for cells moving more slowly closer to a mechanical equilibrium configuration. 
This multirate method is suitable for cell simulations with spatially distributed cell divisions where short time steps are necessary after proliferation but only for the limited number of cells affected by the strong local forces.
In this setting, the coordinates of the dividing cell and its adjacent cells are the fast variables that are integrated with multiple small time steps. The coordinates of the majority of the cells on the other hand are the slow variables and are integrated with a single larger time step chosen as an integer multiple of the small time step such that all cells end up at the same time after the update. 
The dynamic partitioning of the variables into slow and fast is easy in a cell system compared to a general system of ODEs treated in e.g.\ \cite{GeWe84, Logg03}.  Missing values of the slow variables during the short time steps are not interpolated in the multirate method. The method is still of order one \cite{GeWe84}.

While these kinds of methods have not been applied to CBMs before, they have been studied in the context of general systems of ODEs. A systematic variation of the time steps to keep the local discretization errors small for general systems of ODEs is surveyed in \cite{Soderlind06} for single rate methods. 
An early multirate method is found in \cite{GeWe84} where different time steps are taken for different equations in the ODE system.
The variables are partitioned into slow and fast variables with synchronization of the time steps such that the long steps for the slow variables are an integer multiple of the short ones for the fast variables. A conclusion is that the method is competitive for systems with many slow variables.
The global error is bounded adaptively in the multirate method in \cite{Logg03,Logg04} derived from a Galerkin formulation. 
The order of accuracy of the method for different equations may also vary. Multirate time stepping is also developed in \cite{SaHuVe} for systems with different time scales for different sets of equations. The time steps are successively refined by a factor two to satisfy a criterion on the local error. 
For the components integrated by short time steps in a high order method, there will be missing values of the variables advanced by long time steps. These values are obtained by interpolation. 
A recent Runge-Kutta method designed for a fast and a slow time step is found in \cite{SaRoSa19}. 




In this study the following four numerical schemes are compared:
\begin{enumerate}
    \item Single rate forward Euler method
    \item Single rate forward Euler method with stability check
    \item Multirate forward Euler method (with stability check)
    \item Single rate backward Euler method
\end{enumerate}
In order to investigate their properties and the gain from adaptive time stepping in general we apply these schemes to three different cell configurations.
With the same error tolerance for all methods, one method is more efficient than the other ones measured in wall clock time and compared to a method with a constant time step. 
 

The contents of the paper are as follows. The center-based model and the cell forces are given and analyzed in Section~\ref{sec:forceanal}. The numerical methods for time integration of the equations are developed in Section~\ref{sec:adapt_methods} and the conditions for first order accuracy and stability are derived. 
Numerical results for three dimensional cell populations are presented in Section~\ref{sec:numres} and conclusions are drawn in Section~\ref{sec:concl}.

\section{Analysis of forces in center-based models} \label{sec:forceanal}
In a CBM, cells are modelled as intervals, circles, and spheres in one, two and three dimensions (1D, 2D, and 3D). A force on a cell due to the contact with another cell is applied in the center and depends on the distance to the other cell. The system of ODEs defined by these forces governs the motion of the cells. Its properties are investigated in this section. The numerical methods in the next section utilize these properties to obtain accurate and stable time integration of the ODEs.

Matrices and vectors are written in boldface $\fatA$ and $\fatx$. The vector norms are the max norm $\|\cdot\ninf$ and the Euclidean norm $\|\cdot\|_2$ and their associated matrix norms. A time derivative $\d\fatx/\d t$ of $\fatx$ is written $\dot{\fatx}$.

Introduce in $d$ dimensions ($d=1,2,$ or $3$) for $N$ cells the coordinates of the center of cell $i$ $\fatx_i=(x_{i1},\ldots,x_{id})^T,\; i=1,\ldots,N,$  
the relative position $\fatr_{ij}$ to cell $j, j=1,\ldots,N,$ for $j\ne i$, the relative
distance $r_{ij}$, and the direction of the center of the neighbor $\hfatr_{ij}$ 
\begin{equation}
\label{eq:rdef}
      \fatr_{ij}=\fatx_j-\fatx_i,\; r_{ij}=\|\fatr_{ij}\|_2=\left(\sum_{k=1}^d(x_{jk}-x_{ik})^2\right)^{1/2},\; \hfatr_{ij}=r_{ij}^{-1}\fatr_{ij}.
\end{equation}
When $d=3$, $\hfatr_{ij}$ between two cells is defined by two angles $\phi\in[0, 2\pi]$ and $\theta\in[-\pi/2, \pi/2]$ in
\begin{equation}
\label{eq:r3d}
   \hfatr_{ij}=(\cos(\theta)\cos(\phi),\;\cos(\theta)\sin(\phi),\;\sin(\theta))^T.
\end{equation}
For 2D, let $\theta=0$ in \eqref{eq:r3d}. The full coordinate vector for the whole cell system is $\fatx^T=(\fatx_1^T, \fatx_2^T, \ldots,\fatx_N^T)$.

\subsection{A gradient system} \label{sec:gradsyst}

The force on cell $i$ caused by cell $j$ in a CBM is $\hfatr_{ij}g(r_{ij})$ where $\hfatr_{ij}$ is the direction of the force and $\|\hfatr_{ij}g(r_{ij})\|_2=|g(r_{ij})|$ is the modulus of it.
The strength of the force $g(r)$ is continuous and defined for $r\ge 0$ with the properties 
\begin{equation}
\label{eq:gprop}
      0\le r\le s: g(r)\le 0, \quad s<r<r_A: g(r)>0,\quad r\ge r_A: g(r)=0,
\end{equation}
for some positive rest length $s$ and a maximum interaction distance between any two cells $r_A$. The repelling force when $r<s$ is often chosen such that $|g(r)|$ increases when $r$ approaches 0 and $g$ vanishes at $r=s$.
The force is attracting when the distance $r$ is between $s$ and $r_A$.

As an example, take the cubic force defined in \cite{Delile2017MecaGen}
\begin{equation}
g^{\text{cubic}}(r)= \begin{cases}
\mu \left( r - r_A\right)^2 \left( r - s \right) & \text{if } \, 0\le r \leq r_A, \\
0 &\text{otherwise}, \end{cases}
\label{eq:cubic_force}
\end{equation}
where \(\mu\) denotes the spring stiffness. Other examples are found in \cite{MCBH20}.

The potential $G(r)$ for a pair of cells with the force $\hfatr g(r)$ is given by
\begin{equation}
\label{eq:Gdef}
      G(r)=\left\{\begin{array}{ll}-\int_r^s g(\rho)\, d\rho\ge 0,& 0\le r\le s,\\ 
                                   \int_s^r g(\rho)\, d\rho>0,& s< r\le r_A,\\
                                   \int_s^{r_A} g(\rho)\, d\rho=G_A>0,& r> r_A.
                  \end{array}\right.
\end{equation}
The short range potential $G$ is continuously differentiable and $G(s)=0$. The force $\fatf_{ij}$ on cell $i$ due to cell $j$ is
\begin{equation}
\label{eq:Ggrad}
      \fatf_{ij}(\fatx)=-\nabla_{\fatx_i} G(r_{ij})=-\left(\nabla_{\fatx_i} r_{ij}\right) g(r_{ij})=\hat{\fatr}_{ij}g(r_{ij}).
\end{equation}
The simple Hooke's law with $G(r)=\frac{1}{2}(r-s)^2$ and a linear relation between force and distance $g(r)=r-s$ does not satisfy \eqref{eq:Gdef} for $r>r_A$. 
The forces in VBMs are usually defined by the gradient of a potential \cite{VLPJD15}.

The cell system has $N$ free cells and $N_0$ stationary, immobile cells. 
The free cells can move continuously in any direction. The stationary cells do not move in space and have a constant coordinate vector $\fatx_0$. They could be part of a boundary and interact with the free cells via forces as in \eqref{eq:gprop}. 
 
Introduce the potential $V$ including the forces between the free cells, $j=1,\ldots,N,$ and the stationary ones, $j=N+1,\ldots,N+N_0,$
\begin{equation}
\label{eq:potential}
      V(\fatx)=\frac{1}{2} \sum_{i=1}^N \sum_{j=1, j\ne i}^{N+N_0} G(r_{ij}).
\end{equation}
The sum in $V$ has a lower bound 0 due to the properties in \eqref{eq:Gdef}. 
Many terms in the sum are equal to $G_A$ because $r_{ij}>r_A$. 

The contribution of the inertia term in Newton's second law is assumed to be small in a CBM \cite{VLPJD15} but the viscosity term cannot be neglected.
Then the system of $d$ ODEs for the center coordinates $\fatx_i$ of a free cell $i$ is
\begin{equation}
\label{eq:poteqi}
      \dot{\fatx}_i=-\nabla_{\fatx_i} V(\fatx)=\sum_{j=1, j\ne i}^{N+N_0}\fatf_{ij}(\fatx)=\sum_{j=1, j\ne i}^{N+N_0}\hfatr_{ij} g(r_{ij}).
\end{equation}
Since $\hfatr_{ij}=-\hfatr_{ji}$ and $g(r_{ij})=g(r_{ji})$ we have $\fatf_{ij}=-\fatf_{ji}$ following Newton's third law.
There are a limited number of $g(r_{ij})\ne 0$ in \eqref{eq:poteqi} depending on the number of close neighbors to cell $i$ with $r_{ij}<r_A$ that contribute to the total force on the cell.

The ODE system for all free cells is
\begin{equation}
\label{eq:poteq}
      \dot{\fatx}=-\nabla_{\fatx} V(\fatx)=\fatF(\fatx).
\end{equation}
This is a {\it gradient} system with certain properties \cite{Strogatz} and we will see that these properties are inherited by the numerical solution. The force vector $\fatF$ is continuous in $\fatx$ and a continuously differentiable solution $\fatx$ exists and is unique. 

A solution $\fatx_\ast$ is an equilibrium or steady state solution if $\fatF(\fatx_\ast)=0$. An equilibrium solution is not unique for a potential such as \eqref{eq:potential}. It belongs to a subspace $\calD$ with $\fatx_\ast\in\calD$. As an example, consider two free cells, $N=2,\,N_0=0$, with center coordinates $x_1$ and $x_2$ with $x_1<x_2$. Then $\fatx_\ast^T=(x_1, x_1+s)\in\calD$ and $\fatx_\ast^T=(x_1, x_1+\sigma)\in\calD$ for any $x_1$ and any $\sigma>r_A$.

The following theorem characterizes the dynamical system \eqref{eq:poteq} as $t\rightarrow\infty$.

\vspace{2mm}
\noindent
{\bf Theorem 1.}\hspace{2mm} Assume that $V(\fatx)\ge V_{\min}$, $V$  in \eqref{eq:poteq} is continuously differentiable, 
and that $\calD$ is the non-empty set of solutions $\fatx_\ast$ such that $\nabla V(\fatx_\ast)=\fatF(\fatx_\ast)=0$.
Then as $t\rightarrow\infty$ there is no limit cycle and $\fatx(t)\rightarrow \fatx^\infty\in\calD$. 
If $\fatx(0)\in\calD$, then the solution is constant $\fatx(t)=\fatx(0)$ for $t\ge 0$.

\noindent
{\bf Proof.}\hspace{2mm} By \eqref{eq:poteq} we have
\[
      \dot{V}=\dot{\fatx}\cdot\nabla V=-\|\dot{\fatx}\|_2^2=-\nabla V\cdot\nabla V<0,
\]
unless $\fatx\in\calD$. Since $V\ge V_{\min}$, eventually when $t\rightarrow\infty$, $\nabla V(\fatx)\rightarrow 0$, and $\fatx\rightarrow\fatx^\infty\in\calD$.
There is no closed orbit according to \cite[Thm 7.2.1]{Strogatz}. When $\fatx(0)\in\calD$ then $\dot{\fatx}=-\nabla V(\fatx(0))=0$ and the solution is constant.
$\blacksquare$
\vspace{2mm}

\noindent{\bf Remark.}\hspace{2mm} 
The initial condition $\fatx(0)$ determines which $\fatx_\ast$ and constant potential $V_\ast\ge V_{\min}$ the solution converges to. 

\subsection{Frame invariance}
The equations in \eqref{eq:poteq} are frame invariant. To see this, transform the $\fatx_i$ coordinates by a rotation with a $d\times d$ orthonormal matrix $\fatQ$ and a translation $\fatc$ and let $\faty_i=\fatQ\fatx_i+\fatc$. Then define
\[
     \fatrho_{ij}=\faty_j-\faty_i=\fatQ(\fatx_j-\fatx_j),\;
     \rho_{ij}=\|\fatQ(\fatx_j-\fatx_j)\|_2=r_{ij},\; \hfatrho_{ij}=\fatQ\hfatr_{ij}.
\]
The form of the equations is the same after the transformation
\begin{equation}
\label{eq:frameinv}
      \dot{\faty}=\fatQ\dot{\fatx}=\sum_{j=1, j\ne i}^{N+N_0}\fatQ\hfatr_{ij} g(r_{ij})=\sum_{j=1, j\ne i}^{N+N_0}\hfatrho_{ij} g(\rho_{ij}),
\end{equation}
avoiding dependence of the particular coordinate system chosen for the equations. Our numerical methods will also be frame invariant. 

\subsection{Linearization of the forces} \label{sec:linearization}

Next, we consider the linearization of the forces. 
Assume that $V$ is twice continuously differentiable and let $\fatA\in\bbR^{dN\times dN}$ be the Jacobian matrix of $\fatF$ in \eqref{eq:poteq} with elements $A_{k l}=\partial F_k/\partial x_l =-\partial^2 V/\partial x_k\partial x_l,\; k,l=1,\ldots,dN$. The properties of $\fatA$ will be used to ensure stability of our numerical schemes.

After linearization about $\fatx$ for the free cells and $\fatx_0$ for the fixed cells, small perturbations $\dfatx$ of the free cells and $\dfatx_0$ of the fixed cells satisfy
\begin{equation}
\label{eq:dxeq}
      \dot{\dfatx}=\fatA(\fatx)\dfatx+\fatA_F(\fatx,\fatx_0)\dfatx_0.
\end{equation}

The Jacobian matrix $\fatA$ consists of submatrices $\fatA^{ij}\in\bbR^{d\times d},\, i,j=1,\ldots,N,$ for the interaction between the free cells $i$ and $j$. The matrix $\fatA_F\in\bbR^{dN\times dN_0}$ consists of submatrices $\fatA^{ij},\, i=1,\ldots,N,\, j=N+1,\ldots,N+N_0,$ for the interaction between the free cells $i$ and the stationary cells $j$. 
The matrices $\fatA$ and $\fatA_F$ are
\begin{equation}\label{eq:Jacobian}
  \fatA=\left(\begin{array}{ccccc}
  -\sum_{j=1, j\ne 1}^{N+N_0}\fatA^{1j}&\ldots&\fatA^{1j}&\ldots&\fatA^{1N}\\
  &&\vdots&&\\
  \fatA^{i1}&\ldots&-\sum_{j=1, j\ne i}^{N+N_0}\fatA^{ij}&\ldots&\fatA^{iN}\\
  &&\vdots&&\\
  \fatA^{N1}&\ldots&\fatA^{Nj}&\ldots&-\sum_{j=1, j\ne N}^{N+N_0}\fatA^{Nj}
\end{array}\right),
\end{equation}
and
\begin{equation}\label{eq:JacobianF}
  \fatA_F=\left(\begin{array}{ccc}
  \ldots&\fatA^{1j}&\ldots\\
  &\vdots&\\
  \ldots&\fatA^{ij}&\ldots\\
  &\vdots&\\
  \ldots&\fatA^{Nj}&\ldots
\end{array}\right),\; j=N+1,\ldots,N+N_0.
\end{equation}
Each submatrix $\fatA^{ij}$ depends on the force definition as follows
\begin{equation}
\label{eq:subJacobian}
\fatA^{ij}=\hfatr_{ij}\hfatr_{ij}^Tg'(r_{ij})+(\fatI-\hfatr_{ij}\hfatr_{ij}^T)\frac{g(r_{ij})}{r_{ij}},
\end{equation}
where $\fatI$ is the $d\times d$ identity matrix. The Jacobian submatrix depends on the relative positions of the cells and is scaled by the force function and its derivative. 
In \eqref{eq:r3d}, $\hfatr_{ij}$ is represented in 3D by two angles between cells $i$ and $j$. 
Since $\fatA^{ij}$ is symmetric $\fatA^{ij}=\fatA^{ji}=(\fatA^{ij})^T$, $\fatA$ in \eqref{eq:Jacobian} is symmetric with real eigenvalues $\lambda(\fatA)$. 
This is also a consequence of the fact that the equations form a gradient system \eqref{eq:poteq}.

A perturbation of cell $j$ in the direction to the cell $i$, $\dfatx_j=\dx_j\hfatr_{ij}$, induces a correction in the same direction proportional to $g'$
\[
   \fatA^{ij}\dfatx_j=g'(r_{ij})\dx_j\hfatr_{ij}.
\]
With a perturbation in the orthogonal direction, $\dfatx_j=\dx_j(\fatI-\hfatr_{ij}\hfatr_{ij}^T)\fate$ with an arbitrary vector $\fate$ of unit length, the correction in the orthogonal direction is
\[
   \fatA^{ij}\dfatx_j=\frac{g(r_{ij})}{r_{ij}}\dx_j(\fatI-\hfatr_{ij}\hfatr_{ij}^T)\fate.
\] 

If the coordinates $\fatx_j,\, j=N+1,\ldots,N+N_0,$ in $\fatx_0$ are fixed then the perturbation $\dfatx_0=\fatzero$ and only $\fatA$ has an effect on $\dot{\dfatx}$ in \eqref{eq:dxeq}.

If $N_0=0$ there are at least $d$ eigenvectors with eigenvalues 0 of $\fatA$ corresponding to rigid body translation of the system. Let $\fatx^T_\tr=(\tfatx^T,\tfatx^T,\ldots,\tfatx^T)$ with
$d$ different linearly independent $\tfatx\in\bbR^{d}$. Then with $\fatA$ in \eqref{eq:Jacobian} $\fatA\fatx_\tr=\fatzero$ and $\lambda(\fatA)=0$ 
with the associated eigenvector $\fatx_\tr$.
Another transformation of the system without changing the forces between the cells is a rigid body rotation of the whole system or an isolated part of the system as in \eqref{eq:frameinv}.

\subsection{Conservation of the center of gravity of a system of free cells}\label{sec:centgrav}

The center of gravity $\fatx_G$ of a system of free cells is defined by
\begin{equation}
\label{eq:cgrav}
      \fatx_G=\frac{1}{N}\sum_{i=1}^N \fatx_i.
\end{equation}
The behavior of the center of gravity is derived in the next proposition.

\vspace{2mm}
\noindent
{\bf Proposition 1.}\hspace{2mm} 
Assume that the cells in $\calC$ are free and isolated without external forces ($N_0=0$ in \eqref{eq:poteqi}) such that
\begin{equation}
\label{eq:isol}
      \dot{\fatx}_i=\sum_{j\in\calC}\fatf_{ij},\; i\in\calC,\; t\ge 0.
\end{equation}
Then the sum of the coordinates is constant when $t\ge 0$
\begin{equation}
\label{eq:sumcons}
      \sum_{i\in\calC}\sum_{l=1}^d x_{il}(t)=\sum_{i\in\calC}\sum_{l=1}^d x_{il}(0).
\end{equation}
The center of gravity of free cells in \eqref{eq:cgrav} is independent of $t$ 
\begin{equation}
\label{eq:cgravconst}
      \fatx_G(t)=\fatx_G(0).
\end{equation}

\noindent
{\bf Proof.}\hspace{2mm}  The equation \eqref{eq:isol} can be written
\begin{equation}\label{eq:poteqi2}
     \dot{\fatx}=\fatH\fat1_{\calC},
\end{equation}
where $\fatH\in\bbR^{dN\times dN}$ and $\fat1_{\calC}\in\bbR^{dN}$ are defined by
\[
\fatH=\left(\begin{array}{ccccc}
  \fatO&\ldots&\fatF_{1j}&\ldots&\fatF_{1N}\\
  &&\vdots&&\\
  \fatF_{i1}&\ldots&\fatO&\ldots&\fatF_{iN}\\
  &&\vdots&&\\
  \fatF_{N1}&\ldots&\fatF_{Nj}&\ldots&\fatO
\end{array}\right),\quad \fat1_{\calC k}=\left\{\begin{array}{ll}
   1,& k=d(i-1)+l,\\
     & i\in\calC,\; l=1,\ldots,d,\\
   0,& {\rm otherwise}.
\end{array}\right.
\]
The diagonal matrix $\fatF_{ij}\in\bbR^{d\times d}$ consists of the $d$ elements of $\fatf_{ij}$ in the diagonal and $\fatO$ is the $d\times d$ zero matrix.
Since $\fatF_{ij}=-\fatF_{ji}$ according to Newton's third law, $\fatH$ is antisymmetric and $\faty^T\fatH\faty=0$ for any $\faty$. Multiply \eqref{eq:isol} by $\fat1^T_{\calC}$ from the left to
obtain $\fat1^T_{\calC}\dot{\fatx}=\fat1^T_{\calC}\fatH\fat1_{\calC}=0$. Thus, $\frac{d}{dt}\sum_{i\in\calC}\sum_{l=1}^d x_{il}=0$ and \eqref{eq:sumcons} follows.

Let $\fat1_{\calC k}^l=1$ when $k=d(i-1)+l,\; \calC\in\{1,2,\ldots,N\},\; l=1,2,\;{\rm or}\; 3$. Since 
\[
(\fat1_{\calC}^l)^T\dot{\fatx}=\frac{d}{dt}\sum_{i=1}^N x_{il}=(\fat1_{\calC}^l)^T\fatH\fat1_{\calC}^l=0,
\]
we have
\[
     x_{G l}(t)=\frac{1}{N}\sum_{i=1}^Nx_{il}(t)=\frac{1}{N}\sum_{i=1}^Nx_{il}(0)=x_{G l}(0),\quad l=1,\ldots,d.
\]
$\blacksquare$
\vspace{2mm}

An arbitrary sum of coordinates is also preserved by the Euler methods in the next section.


\subsection{Forces after cell division}
\label{sec:forces_after_division}

At a cell proliferation in cell $i$, a new cell is introduced at a given distance from cell $i$. Proliferation at time $t_p$ will make $V$ discontinuous 
\begin{equation}
\label{eq:dpot}
     V(t_p^+)=V(t_p^-)+\Delta V, \; \Delta V>0.
\end{equation}
After $t_p^+$ the cells will move toward an equilibrium as in Theorem 1.
Suppose that the system is at rest at $t_p^-$ such that $r_{ij}=s$ for all $i$ and $j\in\calJ_i$. The set $\calJ_i$ contains the indices of the close neighbors of cell $i$, i.e. if $j\in\calJ_i$ then $r_{ij}<r_A$. 
Let cell $i$ divide into $i$ and $N+1$. Then at $t_p^+$, $\Delta V$ in \eqref{eq:dpot} is 
\begin{equation}
\label{eq:potplus}
      \Delta V=G(r_{i,N+1})+G(r_{N+1,i})+\sum_{j\in\calJ_i}G(r_{ij})+\sum_{j\in\calJ_i}G(r_{N+1,j}).
\end{equation}

After proliferation of the cell at $\fatx_i$ at time $t_p$ as in \eqref{eq:dpot}, a new cell $N+1$ is placed at $\fatx_{N+1}=\fatx_i-\dfatr$ and cell $i$ is moved a short distance to $\fatx_i+\dfatr$ at $t_p^+$.
Then $\fatr_{i,N+1}=-2\dfatr$ and $r_{i,N+1}=2\|\dfatr\|_2=2\dr$. Assume that the cells were in equilibrium before the proliferation with vanishing forces as in \eqref{eq:potplus}. The distances $r_{ij}$ between cell centers are then equal to $s$ in the force model in \eqref{eq:gprop} and \eqref{eq:cubic_force}. 
Furthermore, assume that the distance $2\dr$ between the two new cells is short compared to the distances $r_{ij}$ to the other cells.
Then the forces between cells $i$ and $N+1$ and the other cells in $\calJ_i$ in the sums in \eqref{eq:poteqi} after $t_p$ are to first order in $\dr$ according to \eqref{eq:dxeq} and \eqref{eq:Jacobian}
\begin{equation}
\label{eq:coordplus2}
\begin{array}{rl}
      \dot{\fatx}_i&=\hfatr_{i,N+1} g(2\dr)+\sum_{j\in\calJ_i} \hfatr_{ij} g(r_{ij})+\fatA^{ij}\dfatr\\
                 &=\hfatr_{i,N+1} g(2\dr)-\sum_{j\in\calJ_i} \fatA^{ij}\hfatr_{i,N+1}\dr,\\
      \dot{\fatx}_{N+1}&=-\hfatr_{i,N+1} g(2\dr)+\sum_{j\in\calJ_i} \hfatr_{ij} g(r_{ij})-\fatA^{ij}\dfatr\\
                &=-\hfatr_{i,N+1} g(2\dr)+\sum_{j\in\calJ_i} \fatA^{ij}\hfatr_{i,N+1}\dr.
\end{array}
\end{equation}
Since $|g(2\dr)|$ increases with the force function in \eqref{eq:cubic_force} and the sums in \eqref{eq:coordplus2} decreases when $\dr$ is reduced, the forces in the system are dominated by the force between the two cells involved in the proliferation. This dominance increases the smaller the initial distance $2\dr$ is between the proliferating cells.



\section{Numerical methods for adaptive time integration} \label{sec:adapt_methods}
The ODE system \eqref{eq:poteq} is solved by the forward Euler and  the backward Euler methods for approximations $\fatx^n$ of $\fatx(t^n)$ at the time points $t^{n+1}=t^{n}+\Dt,\; n=0,1,2,\ldots,$ with $t^0=0$. 
The global time step for all cells is $\Dt$ in a single rate algorithm and the initial positions at $t^0$ are $\fatx^0$. The length of the time step is chosen adaptively. Parts of the cell system are integrated by shorter time steps in the multirate algorithm in order to improve the efficiency in particular after proliferation.

The solution is advanced in time by the forward Euler method in
\begin{equation}
\label{eq:forwardE}
      \fatx^{n+1}=\fatx^n+\Dt \fatF(\fatx^n).
\end{equation}
The equation is frame invariant as in \eqref{eq:frameinv}. By multiplying with $\fat1_{\calC}$ as in the proof of the proposition in Section~\ref{sec:centgrav} in a case with $N_0=0$, we find that the coordinate
vector and the center of gravity are conserved by the Euler forward method
\begin{equation}
\label{eq:forwardEcons}
         \sum_{i\in\calC}\sum_{l=1}^d x_{il}^n=\sum_{i\in\calC}\sum_{l=1}^d x_{il}^0,\quad \fatx_G^n=\fatx_G^0,
\end{equation}
cf. \eqref{eq:sumcons}. If $\fatx^n\in\calD$, the space of stationary solutions, then $\fatF(\fatx^n)=0$ and $\fatx^{n+1}=\fatx^n\in\calD$ in \eqref{eq:forwardE}. When $\fatx^n$ has entered $\calD$, it stays there.

\subsection{Single rate time stepping with the forward Euler method} \label{sec:FwInt}

The time steps in this single rate method are controlled such that an accuracy requirement is satisfied in the max norm $\|\cdot\ninf$ and stability is monitored in the Euclidean norm $\|\cdot\|_2$.

\subsubsection{Time step selection} \label{sec:timstsel}

The leading term in the local error in every time step is at $t^n$
\begin{equation}
\label{eq:forwardEerr}
      \fate^{n+1}=\fatx(t^{n+1})-\fatx^{n+1}=\frac{1}{2}\Dt^2\ddot{\fatx}^n=\frac{1}{2}\Dt^2\dot{\fatF}(\fatx^n)=\frac{1}{2}\Dt^2\fatA\fatF(\fatx^n),
\end{equation}
proportional to the acceleration of the cells $\ddot{\fatx}$. When the acceleration vanishes any $\Dt$ will result in a small error.

A numerical approximation of $\fatA\fatF$ in \eqref{eq:forwardEerr} is obtained by one extra evaluation of $\fatF$ in
\begin{equation}
\label{eq:forwardEerrest}
         \fatA\fatF(\fatx^n)\approx\frac{1}{\epsilon}(\fatF(\fatx^n+\epsilon\fatF(\fatx^n))-\fatF(\fatx^n)).
\end{equation}
The parameter $\epsilon$ is small. A discussion how to find it is found in \cite{KnollKeyes}.

Let the error tolerance be $\varepsilon$. Then $\Dt$ is chosen such that the max norm of the local error estimate in each time step in \eqref{eq:forwardEerr} is below this tolerance
\begin{equation}
\label{eq:forwardErelerr}
         \frac{1}{2}\Dt^2\|\fatA\fatF(\fatx^n)\ninf\le\veps.
\end{equation}
The error in each cell coordinate in every time step due to the discretization is then less than $\veps$. The bound on the global time step is then
\begin{equation}
\label{eq:forwardEdt}
         \Dt\le \sqrt{\frac{2\veps}{\|\fatA\fatF\ninf}}.
\end{equation}
With $g$ in \eqref{eq:cubic_force} and by \eqref{eq:poteqi} and \eqref{eq:Jacobian}, we observe that both $\fatA\propto\mu$ and $\fatF\propto\mu$.
Hence, $\Dt\propto\mu^{-1}$ in \eqref{eq:forwardEdt}.

When a limited displacement is allowed in each time step as in \cite{Atwell2016, KBMH, VanLiedekerke2018, SchMH05}, then $\Dt$ satisfies
\begin{equation}
\label{eq:limdisp}
         \|\Dt\dot{\fatx }\ninf=\Dt\|\fatF\ninf\le\veps.
\end{equation}
If $\|\fatF\ninf$ is small then $\Dt$ will be long. This may result in a large error in \eqref{eq:forwardEerrest} if $\|\ddot{\fatx}\ninf$ is not small simultaneously. If $\|\fatF\ninf$ is large but $\|\ddot{\fatx}\ninf$ is small, then unnecessarily short steps are chosen with \eqref{eq:limdisp}.

The procedure is summarized in Algorithm \ref{alg:SRFE}.

\begin{algorithm}[tbp]
\linespread{1.15}\selectfont
\KwIn{Right-hand side \(\textbf{F}\) of Equation \eqref{eq:poteq} defining the ODE system,  start time \(t^0\), final time \(T\), initial coordinates \(\textbf{x}^0\), absolute accuracy \(\varepsilon\), approximation parameter for Jacobian-force product \(\epsilon\)}
Initialize \(t = t^0\); \(\textbf{x} = \textbf{x}^0\)\;
\While{\(t<T\)}{
Evaluate force for current coordinates \(\hat{\textbf{F}} = \textbf{F}(\textbf{x})\)\;
Approximate Jacobian-force product by \(\textbf{AF} = \frac{1}{\epsilon} \left(\textbf{F}( \textbf{x}+ \epsilon\, \hat{\textbf{F}})-\hat{\textbf{F}}\right)\)\;
Calculate time step size \(\Delta t=\sqrt{\frac{2\,\varepsilon}{ \Vert \textbf{AF} \ninf}}\)\;
Update \(\textbf{x} \leftarrow \textbf{x} + \Delta t \, \hat{\textbf{F}}\); \(t \leftarrow t + \Delta t\)\;
}
\SetAlgoRefName{I}
\caption{Single rate forward Euler method (SRFE)}
\label{alg:SRFE}
\end{algorithm}

\subsubsection{Explicit consideration of the stability bound}\label{sec:explstab}
The difference between $\fatx^n$ and $\fatx^{n+1}$ fulfills
\begin{equation}
\label{eq:forwardEdiff}
\begin{array}{rll}
  \fatx^{n+1}-\fatx^n&=\fatx^{n}-\fatx^{n-1}+\Dt(\fatF(\fatx^{n})-\fatF(\fatx^{n-1}))\\
                    &=\fatx^{n}-\fatx^{n-1}+\Dt\tfatA(\fatx^{n}-\fatx^{n-1})\\
                    &=(\fatI+\Dt\tfatA)(\fatx^{n}-\fatx^{n-1}).
\end{array}
\end{equation}
The matrix $\tfatA$ depends on $\fatx^{n-1}, \fatx^n,$ and $\fattheta$ with components $\theta_k,\; k=1,\ldots, dN,$ between 0 and 1. 
By the mean value theorem, each row $k$ of $\tfatA$ is evaluated as $\nabla_\fatx F_k(\theta_k\fatx^{n-1}+(1-\theta_k)\fatx^n)$.

The symmetric $\fatA$ in \eqref{eq:Jacobian} has a factorization
\begin{equation}
\label{eq:Afact}
      \fatA=\fatQ\fatLambda\fatQ^T=(\fatQ_+, \fatQ_-)\left(\begin{array}{cc}\fatLambda_+&0\\0&\fatLambda_-\end{array}\right)(\fatQ_+, \fatQ_-)^T.
\end{equation}
The positive eigenvalues of $\fatA$ are in the diagonal of $\fatLambda_+$ and the non-positive eigenvalues in the diagonal of $\fatLambda_-$. The square, orthonormal eigenvector matrix $\fatQ$ is partitioned into $(\fatQ_+, \fatQ_-)$ corresponding to the partitioning of the diagonal eigenvalue matrix $\fatLambda$. 

In the difference between successive solutions in \eqref{eq:forwardEdiff},  $\tfatA$ is approximated by $\fatA$ at $\fatx^n$. If $\Dt$ is small then the difference $\tfatA-\fatA$ is small. For the components corresponding to non-positive eigenvalues, we have
\begin{equation}
\label{eq:succdiffnpos}
\begin{array}{rl}      \fatQ_-^T(\fatx^{n+1}-\fatx^n)&=\fatQ_-^T(\fatI+\Dt\fatA)(\fatx^{n}-\fatx^{n-1})\\
&=(\fatI+\Dt\fatLambda_-)\fatQ_-^T(\fatx^{n}-\fatx^{n-1}).
\end{array}
\end{equation}
A stability requirement on the numerical integration is that 
\begin{equation}
\label{eq:stabreq}
\|\fatQ_-^T(\fatx^{n+1}-\fatx^n)\|_2\le \|\fatQ_-^T(\fatx^{n}-\fatx^{n-1})\|_2
\end{equation}
in the Euclidean norm. 
A sufficient condition for stability in \eqref{eq:stabreq} is then that $\Dt$ is such that
\begin{equation}
\label{eq:forwardEstab}
      \|\fatI+\Dt\fatLambda_-\|_2\le 1.
\end{equation}
The equivalent condition on $\Dt$ is that for all non-positive eigenvalues $\lambda_{k-}$ of $\fatA$ in $\fatLambda_-$
\begin{equation}
\label{eq:forwardEstab2}
      -2\le \Dt\lambda_{k-}\le 0.
\end{equation}
Either $\Dt$ is restricted by the stability constraint in $\ell_2$ \eqref{eq:forwardEstab2} or the accuracy constraint in $\ell_\infty$ \eqref{eq:forwardEdt}.

The eigenvalues of $\fatA$ in the stability bound in \eqref{eq:forwardEstab2} can be estimated by Gershgorin's theorem \cite[7.2.1]{Lanc}. The eigenvalues are located in the union of the intervals  $[\xi_k-\rho_k, \xi_k+\rho_k],\, k=1,\ldots ,dN,$ 
where
\begin{equation}
\label{eq:Gersh}
   \xi_k=A_{kk}=-\sum_{j=1, j\ne i}^N A^{ij}_{ll},\quad \rho_k=\sum_{m=1, m\ne k}^{dN}|A_{km}|,
\end{equation}
according to \eqref{eq:Jacobian}. The indices in \eqref{eq:Gersh} are related by $k=d(i-1)+l,\; i=1,\ldots,N, \; l=1,\ldots,d$.

The leftmost eigenvalue $\lambda_L=\min_k \lambda_k(\fatA)<0$ constrains $\Dt$ in \eqref{eq:forwardEstab2} such that $\Dt\le 2/|\lambda_L|$ and
can be estimated by $\min_k \xi_k-\rho_k\le\lambda_L(\fatA)$ with a tighter bound $2/|\min_k \xi_k-\rho_k|\le 2/|\lambda_L|$ on $\Dt$.
Some eigenvalues lie on the positive real axis corresponding to growing modes.
There the accuracy puts a bound on $\Dt$ such that the local error in each step in \eqref{eq:forwardEerr}
is less than some $\varepsilon$.

Gershgorin's estimate in \eqref{eq:Gersh} depends on the number of non-zero terms $A_{km}$. That number depends on the number of neighboring cells in $\calJ_i$. Its maximum $N_d$ in a cell system depends on $d$ but is independent of $N$ when $N$ is large. 
In a crowded aggregation in 2D, the maximum number of circles of equal size touching a circle in the middle is six (the kissing number \cite{Conway}). At most twelve surrounding spheres touch a sphere in the middle in 3D.
The elements $A_{km}$ depend on the direction angles of $\hfatr$, $g'(r)$, and $g(r)/r$.
In general, the Gershgorin estimate of the minimal eigenvalue  $\min_k\xi_k-\rho_k$ decreases when the dimension increases because of an increasing number of non-zero terms in the sums. 

The eigenvalue of largest modulus of the Jacobian matrix $\fatA$ is estimated in \cite{HaWa, Petzold83} in numerical methods for ODEs where several, successive evaluations of $\fatF$ are available in each step mimicking a power iteration. With the forward Euler method we have only a few evaluations at our disposal. 

Algorithm \ref{alg:SRFES} is the integration method with a global time step adapted to the solution as in Algorithm \ref{alg:SRFE} but with a stability constraint.   

For certain geometries of the cell system, the eigenvalues are available explicitly. Consider a Cartesian configuration of cells in dD with $N$ cells in each coordinate direction with $\hfatr_{ij}$ pointing along the coordinate axes.
The system is close to steady state.
Then in 3D we have $\theta=-\frac{\pi}{2}, 0, \frac{\pi}{2}$, and for $\theta=0$ we let $\phi=0, \frac{\pi}{2}, \pi, \frac{3\pi}{2}$ in \eqref{eq:r3d} and in 2D $\theta=0$ with the same $\phi$ as in 3D. The outer layer of cells has fixed positions.
Linearize around the equilibrium $r_{ij}=s$ to obtain $\fatA$.
The coordinates in $\fatx$ are ordered such that all $x$ coordinates are first, then all $y$ coordinates come, and finally in 3D the $z$ coordinates.
Then $\fatA$ consists of $d$ tridiagonal submatrices $g'(s)\fatB$ on the diagonal, one for each coordinate. Each $\fatB$ has $N^{(d-1)}$ tridiagonal matrices $\fatB_0$ on
the diagonal where $\fatB_0$ has $-2$ on the diagonal and $1$ on the subdiagonal and the superdiagonal.
The eigenvalues of $\fatB_0$ (and $\fatB$) are $\lambda_j=-2(1-\cos(j\pi/(N+1)))<0,\; j=1,2,\ldots,N$. The eigenvalue of $\fatA$ of largest modulus is therefore for large $N$
\begin{equation}
\label{eq:eiglarge}
    \lambda_N=\left(-2+2\cos(\frac{N\pi}{N+1})\right)g'(s)\approx \left(-4+\frac{\pi^2}{(N+1)^2}\right)g'(s).  
\end{equation}
The estimate with Gershgorin's theorem is $\lambda_N=-4g'(s)$ which is the limit in \eqref{eq:eiglarge}  for large $N$. The time step restriction for stability in \eqref{eq:forwardEstab2} is then $\Dt\le 1/(2g'(s))$. 
The eigenvalue of smallest modulus and closest to zero is 
\begin{equation}
\label{eq:eigsmall}
      \lambda_1=-2\left(1-\cos(\frac{\pi}{N+1})\right)g'(s)\approx -\frac{\pi^2}{(N+1)^2}g'(s).
\end{equation}

\begin{algorithm}[tbp]
\linespread{1.15}\selectfont
\KwIn{Right-hand side \(\textbf{F}\) of Equation \eqref{eq:poteq} defining the ODE system,  start time \(t^0\), final time \(T\), initial coordinates \(\textbf{x}^0\), absolute accuracy \(\varepsilon\), Jacobian \(\textbf{A}\)}
Initialize \(t = t^0\); \(\textbf{x} = \textbf{x}^0\)\;
\While{\(t<T\)}{
Evaluate force and Jacobian for current coordinates \(\hat{\textbf{F}} = \textbf{F}(\textbf{x})\); \(\hat{\textbf{A}} = \textbf{A}(\textbf{x})\)\;
Estimate smallest eigenvalue \(\lambda_{\text{min}}\) of \(\hat{\textbf{A}}\) using Gershgorin's theorem\;
Calculate time step size \(\Delta t= \min \left(\sqrt{\tfrac{2\,\varepsilon}{\Vert \hat{\textbf{A}}\hat{\textbf{F}} \ninf }}, \tfrac{2}{|\lambda_{\text{min}}|} \right)\)\;
Update \(\textbf{x} \leftarrow \textbf{x} + \Delta t \, \hat{\textbf{F}}\); \(t \leftarrow t + \Delta t\)\;
}
\SetAlgoRefName{II}
\caption{Single rate forward Euler method with stability bound (SRFES)}
\label{alg:SRFES}
\end{algorithm}

\subsubsection{Restriction on the time step size after cell proliferation}
\label{sec:restr_step_size_proliferation}

The linearization matrix $\fatA$ after proliferation as in  \eqref{eq:coordplus2} consists of a matrix $\fatA_0$ with contribution from all forces except for the force between the divided cell $i$ and the new cell $N+1$ and a matrix $\fatA_p$ such that $\fatA=\fatA_0+\fatA_p$. 
Then the non-zero submatrices of the proliferation matrix $\fatA_p$ are
\begin{equation}
\label{eq:Aprolifeig}
      \left(\begin{array}{cc}-\fatA^{i,N+1}&\fatA^{i,N+1}\\ \fatA^{ N+1,i}&-\fatA^{N+1,i}\end{array}\right).
\end{equation}

At cell proliferation in \eqref{eq:poteqi} with $\fatn=\hfatr_{i,N+1}$ and $r_0=2\dr<s$ in \eqref{eq:coordplus2}, $\fatF_i\approx \fatn g(r_0)$ and $\fatF_{N+1}\approx -\fatn g( r_0)$ when $\dr$ is small. Moreover,
\begin{equation}
\label{eq:Aprolif}
    \fatA^{i,N+1}=\fatA^{N+1,i}=\fatn\fatn^T g'(r_0)+(\fatI-\fatn\fatn^T)g(r_0)/ r_0
\end{equation}
in \eqref{eq:subJacobian}. Thus,
\begin{equation}
\label{eq:AFprolif}
         (\fatA\fatF)_i\approx -2\fatn g'( r_0)g( r_0),\; (\fatA\fatF)_{N+1}\approx 2\fatn g'( r_0)g( r_0)
\end{equation}
and 
$\|\fatA\fatF\ninf\approx 2g'( r_0)|g( r_0)|\max_l|n_l|$.  

In general, the eigenvalues of largest modulus immediately after proliferation of $\fatA_p$ are much larger than the modulus of the eigenvalues of $\fatA_0$, see \eqref{eq:coordplus2}. It follows from \cite[Ch. 8.1]{GVL} that for $k=1,\ldots,dN$, 
\begin{equation}
\label{eq:Aprolifeigest}
      \lambda_k(\fatA_p)+\min_j(\lambda_j(\fatA_0))\le\lambda_k(\fatA)\le\lambda_k(\fatA_p)+\max_j(\lambda_j(\fatA_0)).
\end{equation}
Since $|\lambda(\fatA_p)|\gg|\lambda(\fatA_0)|$ in \eqref{eq:Aprolifeigest} the dominant eigenvalues of $\fatA$ after cell proliferation in \eqref{eq:Aprolif} are close to those in $\fatA_p$.

When $d=3$ and a force similar to \eqref{eq:cubic_force} with $g'(r_0)>0$ for $r_0<s$, the eigenvalues and eigenvectors of the submatrix of $\fatA_p$ in \eqref{eq:Aprolifeig} are by \eqref{eq:subJacobian}
\begin{equation}
\label{eq:Aprolifeig2}
\begin{array}{lll}
   \lambda_1=0,\; (\fatn^T, \fatn^T),& \lambda_2=-2g'( r_0)<0,\; (\fatn^T, -\fatn^T),\\  
   \lambda_3=0,\; (\fatm_1^T, \fatm_1^T),& \lambda_4=-2g( r_0)/ r_0>0,\; (\fatm_1^T, -\fatm_1^T),\quad \fatm_1^T\fatn=0,\\  
   \lambda_5=0,\; (\fatm_2^T, \fatm_2^T),& \lambda_6=-2g( r_0)/ r_0>0,\; (\fatm_2^T, -\fatm_2^T),\\ &\fatm_2^T\fatn=0,\; \fatm_1^T\fatm_2=0.  
\end{array}
\end{equation}
The other $3N-6$ eigenvalues of $\fatA_p$ are zero.
The interpretation of the second eigenvalue and eigenvector is that a small perturbation of the position of cell $i$ in the $\fatn$ direction from cell $i$ to cell $N+1$ increases the force on cell $i$ in the opposite direction by $2g'(r_0)$ in a repellation. 
A perturbation in the plane orthogonal to $\fatn$ increases the force in the same direction with the strength of eigenvalues four and six.
The first four eigenvalues and eigenvectors in 2D are the same as in \eqref{eq:Aprolifeig2}.

The restriction on the time step for stability for negative eigenvalues in $\|\cdot\|_2$ and accuracy in $\|\cdot\ninf$ after proliferation is then by \eqref{eq:forwardEdt}, \eqref{eq:AFprolif}, \eqref{eq:forwardEstab2}, and \eqref{eq:Aprolifeig2} 
\begin{equation}
\label{eq:Aprolifrestr}
      \Dt\le\min\left(\frac{1}{g'( r_0)},\left(\frac{\varepsilon}{g'( r_0)|g( r_0)|\max|n_l|}\right)^{1/2}\right).
\end{equation}
If $\varepsilon g'(r_0)/|g( r_0)|\max|n_l|<1$, then after proliferation $\Dt$ in \eqref{eq:Aprolifrestr} is determined by the accuracy constraint and is independent of $N$ and the dimension. 



\subsection{Multirate adaptive time stepping with the forward Euler method} \label{sec:localerr}

The time steps are chosen differently in different parts of the cell system in this multirate method. They satisfy an accuracy bound \eqref{eq:forwardEdt} locally with small steps close to a proliferation and larger steps in quiescent parts of the system. In addition, they also satisfy the stability bound in \eqref{eq:forwardEstab2}.


\subsubsection{Time step selection}\label{sec:MRFE}

Let the acceleration in the leading term in the local error in \eqref{eq:forwardEerr} be denoted by $\fateta=\fatA\fatF$ and let $\calK$ be the set of equations $\calK=\{1,2,\ldots,dN\}$. 
Introduce the two subsets $\calK_\kappa, \, \kappa=0,1,$ and the time steps $\Dtau_\kappa$ for all equations in each $\calK_\kappa$. The method is generally applicable but is suitable for ODE systems where the error estimate requires small time steps occasionally and for a limited number of equations. The sets $\calK_\kappa$ are disjunct and cover $\calK$, $\calK=\calK_0\bigcup\calK_1$. The relation between the time steps is chosen to be 
\begin{equation}
\label{eq:taudef} 
    \Dtau_1=m\Dtau_0,
\end{equation}
with an integer $m> 1$ such that $\Dtau_0<\Dtau_1$. The coordinate of a cell with an equation $k$ in $\calK_\kappa$ is advanced in time by $\Dtau_\kappa$. After $m$ steps with $\Dtau_0$ in $\calK_0$ and one step with $\Dtau_1$ in $\calK_1$ from $t^n$, all cell coordinates have reached the same $t^{n+1}=t^n+\Dtau_1$. Let $\fateta_\kappa=\{\eta_k|k\in\calK_\kappa\}$. Then the local errors $\fate_\kappa$ in $\calK_\kappa$ in this step are by \eqref{eq:forwardEerr} approximately
\begin{equation}
\label{eq:adaplocerr} 
    m\fate_0=\frac{m}{2}\fateta_0\Dtau_0^2
            =\frac{1}{2m}\fateta_0\Dtau_1^2,\quad
    \fate_1=\frac{1}{2}\fateta_1\Dtau_1^2.
\end{equation}
The time steps are chosen to satisfy an accuracy bound $\veps$ on the local error  
\begin{equation}
\label{eq:locerreps} 
    \|m\fate_0\ninf\le\veps,\quad
    \|\fate_1\ninf\le\veps.
\end{equation}
The time steps are either bounded by the accuracy of the forward Euler method as in \eqref{eq:locerreps} or the stability of the method as in \eqref{eq:forwardEstab2}.

Introduce $\chi_{a0}, \chi_{a1}, \chi_{s0},$ and $\chi_{s1}$.
The time steps $\Dtau_0$ and $\Dtau_1$ are determined by the parameters $\chi_{a\kappa}$ if they are bounded by the accuracy and by $\chi_{s\kappa}$ if bounded by stability.  

The maximum of $|\eta_k|$ in $\calK$ defines $\chi_{a0}$ and $\chi_{a1}$
\begin{equation}
\label{eq:chidefacc}
      \chi_{a0}=\max_{k\in\calK}|\eta_k|,\;\chi_{a1}=\frac{\chi_{a0}}{m}.
\end{equation}
The maximum time step $\Dtaus$ for stability with the Euler forward method is given by \eqref{eq:forwardEstab} and \eqref{eq:forwardEstab2}: $\Dtaus=2/|\min_k(\lambda_k(\fatA)|$. Define $\chi_{s0}$ and $\chi_{s1}$ using $\Dtaus$ 
\begin{equation}
\label{eq:chidefstab}
      \chi_{s1}=\frac{2\veps}{\Dtaus^2}=\frac{\veps}{2}(\min_k(\lambda_k(\fatA))^2, \;
      \chi_{s0}=m\chi_{s1}.
\end{equation}

If the time steps are constrained by the accuracy, $\chi_{a1}\ge\chi_{s1}$, then let $\chi_{1}=\chi_{a1}$ and if the stability constrains the time steps, $\chi_{a1}<\chi_{s1}$, then let $\chi_1=\chi_{s1}$. Introduce $\chi_0=m\chi_1$ and the sets $\calK_0$ and $\calK_1$ as follows 
\begin{equation}
\label{eq:Kdef1}
      \calK_0=\{k| \;\chi_1<|\eta_k|\le\chi_0\},\;
      \calK_1=\{k| \;|\eta_k|\le\chi_1\},
\end{equation}
and let
\begin{equation}
\label{eq:taudef1}
      \Dtau_0=\sqrt{\frac{2\veps}{m\chi_{0}}},\;
      \Dtau_1=\sqrt{\frac{2\veps}{\chi_{1}}}.
\end{equation}

Firstly, consider the case when $\chi_{a1}\ge\chi_{s1}$. The time step is bounded by the accuracy. If $\calK_1\ne\emptyset$ then there is at least one $k$ with $|\eta_k|\le\chi_{a1}$. Choose the time steps in $\calK_0$ and $\calK_1$ such that
\begin{equation}
\label{eq:tau01def}
      \Dtau_0=\sqrt{\frac{2 \veps}{m\chi_{a0}}},\;
      \Dtau_1=\sqrt{\frac{2 \veps}{\chi_{a1}}}=m\Dtau_0=\sqrt{\frac{2 m \veps}{\chi_{a0}}}.
\end{equation}
Since $\chi_{s1}\le\chi_{a1}<\chi_{a0}$ it follows from \eqref{eq:adaplocerr} that both the accuracy and the stability requirements $\Dtau_0<\Dtau_s$ and $\Dtau_1\le\Dtau_s$ in \eqref{eq:locerreps} and \eqref{eq:forwardEstab2} are fulfilled. 
If $\calK_1=\emptyset$ then all cells are advanced by $\Dtau_0$.

Secondly, assume that $\chi_{a1}<\chi_{s1}$. Then stability bounds the time step. If for all $k$ we have $\chi_{s1}<|\eta_k|\le\chi_{a0}$ then $\calK_1=\emptyset$. The time step is then $\Dtau_0=\sqrt{2\veps/m\chi_{s0}}<\Dtau_s$. If $\calK_1\ne\emptyset$ and $\calK_0\ne\emptyset$ then $\chi_{a1}<\chi_{s1}\le\chi_{a0}$ and $\Dtau_0$ and $\Dtau_1$ are
\begin{equation}
\label{eq:tau02def}
      \Dtau_0=\sqrt{\frac{2 \veps}{m\chi_{s0}}}\le\sqrt{\frac{2 \veps}{m\chi_{a0}}},\;
      \Dtau_1=\sqrt{\frac{2 \veps}{\chi_{s1}}}=\Dtau_s<\sqrt{\frac{2  \veps}{\chi_{a1}}}.
\end{equation}
It follows from \eqref{eq:tau02def} that both $\Dtau_0$ and $\Dtau_1$ satisfy the accuracy requirements. In case $\chi_{a0}<\chi_{s1}$ then $\calK_0=\emptyset$ and all cells are integrated by $\Dtau_1=\sqrt{2\veps/\chi_{s1}}=\Dtau_s$.

The solution is advanced from $t^{n}$ to $t^{n+1}$ by the multirate forward Euler method (MRFE) in Algorithm \ref{alg:MRFE}. After each small time step \(\Dtau_0\) the force terms \(\hat{\fatF}^{(0)}\) for equations in \(\calK_0\) need to be updated. After \(m\) small time steps the entries in  \(\hat{\fatF}^{(1)}\) that have been affected by an update to \(\fatx^{(0)}\) need to be updated as well before updating \(\fatx^{(1)}\). 

The global order of accuracy of the time integration is one because the global errors in the values in $\calK_1$  needed in $\calK_0$ are of $\calO(\Dtau_1)$ and those needed in $\calK_1$ are also of $\calO(\Dtau_1)$, see \cite[p. 490]{GeWe84}. 
If all equations are integrated with the same time step, $\calK_0=\emptyset$ or $\calK_1=\emptyset$, then the method is the usual single rate forward Euler method.

\begin{algorithm}[tbp]
\linespread{1.15}\selectfont
\KwData{Right-hand side \(\textbf{F}\) of Equation \eqref{eq:poteq} defining the ODE system,  start time \(t^0\), final time \(T\), initial coordinates \(\textbf{x}^0\), absolute accuracy \(\varepsilon\), Jacobian \(\textbf{A}\), ratio between levels \(m\)}
Initialize \(t = t^0\); \(\textbf{x} = \textbf{x}^0\)\;
\While{\(t<T\)}{
Evaluate force and Jacobian for current coordinates \(\hat{\textbf{F}} = \textbf{F}(\textbf{x})\); \(\hat{\textbf{A}} = \textbf{A}(\textbf{x})\)\;
Estimate smallest eigenvalue \(\lambda_{\text{min}}\) of \(\hat{\textbf{A}}\) using Gershgorin's theorem\;
Calculate time step sizes \(\Delta \tau_1= \min \left(\sqrt{\tfrac{2\,\varepsilon m}{\Vert \hat{\textbf{A}}\hat{\textbf{F}} \Vert_\infty }}, \tfrac{2}{|\lambda_{\text{min}}|} \right)\); \(\Delta \tau_0 = \Delta \tau_1/m\)\;
Split \(\textbf{x}\) into \(\textbf{x}^{(0)}\) and \(\textbf{x}^{(1)}\) and $\hat{\fatF}$ into $\hat{\fatF}^{(0)}$ and $\hat{\fatF}^{(1)}$ using $\calK_0$ and $\calK_1$ taking \(\chi_1 = \frac{2\,\varepsilon}{(\Delta \tau_1)^2}\) in \eqref{eq:Kdef1}\; 
\For{\(j = 0; j<m; j++\)}{
Update \(\textbf{x}^{(0)} \leftarrow \textbf{x}^{(0)} + \Delta \tau_0 \, \hat{\textbf{F}}^{(0)}\)\;
Update \( \hat{\textbf{F}}^{(0)} \leftarrow \textbf{F}(\textbf{x}^{(0)})\) \; 
}
Do partial update \( \hat{\textbf{F}}^{(1)} \leftarrow \textbf{F}(\textbf{x}^{(0)}, \textbf{x}^{(1)})\) for all entries affected by changes to \(\textbf{x}^{(0)}\)\;
Update \(\textbf{x}^{(1)} \leftarrow \textbf{x}^{(1)} + \Delta \tau_1 \, \hat{\textbf{F}}^{(1)}\)\;
Assemble \(\textbf{x}\) from \(\textbf{x}^{(0)}\) and \(\textbf{x}^{(1)}\)\;
Update \(t \leftarrow t + \Delta \tau_1\)\;
}
\SetAlgoRefName{III}
\caption{Multirate forward Euler method (MRFE)}
\label{alg:MRFE}
\end{algorithm}

\subsubsection{Convergence of adaptive scheme}

Between $t^n$ and $t^{n+1}$, the local errors in each step with $\Dtau_0$ or $\Dtau_1$ are approximately summed as assumed in \eqref{eq:adaplocerr} and shown in the theorem below.
It follows from a proposition that the stability condition in \eqref{eq:forwardEstab2} is satisfied also by the separate steps in the multirate method.

At the fine level, $m$ steps of length $\Dtau_0$ are taken and at the coarse level one step of length $\Dtau_1$.
Since $\max_k|\eta_k|\le \chi_\kappa$ in $\calK_\kappa$ in \eqref{eq:Kdef1}, the maximum local errors $\fate_\kappa$ in $\calK_\kappa$ at the new time level $t^{n+1}=t^n+\Dtau_1$ are
\begin{equation}
\label{eq:localseterr}
     \|\fate_1\ninf=\frac{1}{2}\Dtau_1^2\chi_{a1},\quad
      \|\fate_0\ninf\approx\frac{m}{2}\Dtau_0^2\chi_{a0}.
\end{equation}

The next theorem is a modification of a theorem in \cite{HaNoWa} for an explicit integration method with a global time step $\Dt$ to compute $\fatx^{n+1}$ from $\fatx^n$ with an increment function $\fatPhi$ as in
\begin{equation}
\label{eq:explRK}
    \fatx^{n+1}=\fatx^{n}+\Dt\fatPhi(t^n, \fatx^n, \Dt). 
\end{equation}
It is applicable to any scheme of this form, e.g. an explicit Runge-Kutta method.

\vspace{2mm}
\noindent
{\bf Theorem 2.}\hspace{2mm} For the integration method in \eqref{eq:explRK} with the Jacobian $\fatA$ of $\fatF$ in \eqref{eq:poteq}, assume that in a neighborhood of $\fatx(t),\, t\in[0, T]$, there is a bound 
\begin{equation}
\label{eq:Ldef}
    \max_{j=1,\ldots,dN} \,(A_{jj}+\sum_{k=1, k\ne j}^{dN} |A_{jk}|)\le L,
\end{equation}
and that the variable time steps $\Dt_i=t^{i+1}-t^i,\, i=0,1,\ldots,n-1,$ are so small that the numerical solution stays in this neighborhood. The local error is bounded by $\|\fate^i\ninf\le C_i\Dt_i^{p+1}$ where $0<C_{\min}\le C_i\le C_{\max}$ and $0<\Dt_i\le\Dt_{\max}$. The time steps are chosen such that $C_i\Dt_i^{p+1}=\veps$. The global error $\fatE^n$ at $t^n$ is defined by the exact solution $\fatx(t)$ and the numerical solution $\fatx$ in $\fatE^n=\fatx(t^n)-\fatx^n$. The global error satisfies
\[
     \|\fatE^n\ninf\le \varepsilon^{p/(p+1)}\frac{C_{\max}}{C_{\min}^{p/(p+1)}}\frac{C'}{L}\left(\exp(Lt^n)-1\right),
\]
where $C'=1$ if $L\ge 0$ and $C'=\exp(-L\Dt_{\max})$ if $L<0$.

\noindent
{\bf Proof.}\hspace{2mm} 
The bound in \eqref{eq:Ldef} is a bound on the logarithmic norm of $A$ in the maximum norm, see \cite[Th. I.10.5]{HaNoWa}. The global error in Theorem II.3.4 in \cite{HaNoWa} is bounded by 
\begin{equation}
\label{eq:Ebound}
    \|\fatE^n\ninf\le \Dt_{\max}^p C_{\max}\frac{C'}{L}\left(\exp(Lt^n)-1\right).
\end{equation}
Since $\Dt_{\max}=\left(\varepsilon/C_{\min}\right)^{1/(p+1)}$ in \eqref{eq:Ebound}, the estimate in the theorem follows.
$\blacksquare$

{\bf Remark.}\hspace{2mm} The bound in \eqref{eq:Ldef} is the same as in the Gershgorin estimate of the maximum eigenvalue in \eqref{eq:Gersh} $\xi_k+\rho_k$. The increment function in the Euler forward method  \eqref{eq:forwardE} is $\fatPhi=\fatF(\fatx^n)$ in \eqref{eq:explRK}. 

Apply the theorem with $\Dt=\sqrt{\varepsilon/C}=\Dtau_0, p=1, C=\frac{1}{2}\chi_{a0}, t^1=m\Dtau_0=\Dtau_1,$ and $L$ as in \eqref{eq:Ldef} to determine the local error after $m$ steps on the fine level. In most of the interval $[0, T]$ between proliferations, $L$ is small and positive. If $L$ is small and $m$ not too large, then by \eqref{eq:Ebound}
\[
  \|\fatE^1\ninf\le\frac{\chi_{a0}\Dtau_0}{2L}(\exp(L m\Dtau_0)-1)\approx\frac{1}{2}m\chi_{a0}\Dtau_0^2=m\|\fate_0\ninf.
\]
Hence, the assumption in \eqref{eq:adaplocerr} on successive short steps $\Dtau_0$ holds true.

If the time step is bounded by accuracy then by Theorem 2
\begin{equation}
\label{eq:localseterr2}
\begin{array}{rl}           
      \|\fate_1\ninf&\approx\displaystyle{\frac{1}{2}\Dtau_1^2\chi_{a1}=\frac{m^2}{2}\Dtau_0^2\chi_{a1}=\frac{m}{2}\Dtau_0^2\chi_{a0}= \veps,}\\
      m\|\fate_0\ninf&\approx\displaystyle{\frac{m}{2}\Dtau_0^2\chi_{a0}=\veps}.
\end{array}
\end{equation}
The local error $\|\fate\ninf=\max(m\|\fate_0\ninf, \|\fate_1\ninf)$ in the multirate method in Section \ref{sec:localerr} fulfills the error criterion $\|\fate\ninf\le\veps$ between $t^n$ and $t^{n+1}$. The theorem can be applied to the locally adaptive method with $\Dt=\Dtau_1$ and a given error tolerance $\veps$. Numerical examples in Section \ref{sec:numres} with Algorithm \ref{alg:MRFE} and $p=1$ confirm the dependence of $\veps$ in $\|\fatE^n\ninf$ in the theorem. Theorem 2 is also directly applicable to the solutions obtained with Algorithms \ref{alg:SRFE} and \ref{alg:SRFES} with the same $\Dt$ for all equations.



The next proposition shows that if the eigenvalue condition \eqref{eq:forwardEstab2} is satisfied by $\fatA$ as assumed to obtain the stable time step, then the eigenvalue condition is also satisfied by the time steps $\Dtau_0$ and $\Dtau_1$. 

The difference between two consecutive $\fatx^n$ in Algorithm \ref{alg:MRFE} is propagated from $t^n$ to $t^{n+1}=t^n+\Dtau_1$ by first computing the variables in $\calK_0$ and then using these values to update the variables in $\calK_1$
\begin{equation}
\label{eq:multistab}
   \fatx^{n+1}-\fatx^n=(\fatI+\Dtau_1\fatA_1)(\fatI+\Dtau_0\fatA_0)^{m}(\fatx^{n}-\fatx^{n-1}).
\end{equation}
The matrix $\fatA_\kappa$ consists of the rows of $\fatA$ with an index $k$ in $\calK_\kappa$. The remaining rows in $\fatA_\kappa$ are zero. The number of indices in $\calK_\kappa$ is denoted by $N_\kappa$. The $N_\kappa$ non-zero eigenvalues of $\fatA_\kappa$ are bounded from below and above by the eigenvalues of $\fatA$ in the next proposition.

{\bf Proposition 2.} 
The eigenvalues of the symmetric $\fatA$ are in the interval $\calI_\fatA=[\lambda_{\fatA \min}, \lambda_{\fatA \max}]$.
There are $N-N_\kappa$ zero eigenvalues of $\fatA_\kappa$. 
The non-zero eigenvalues of $\fatA_\kappa$ are in the interval $\calI_{\fatA_\kappa}=[\lambda_{{\fatA_\kappa} \min}, \lambda_{{\fatA_\kappa} \max}]$. 
Then $\calI_{\fatA_\kappa}\subseteq\calI_\fatA$, i.e. 
\begin{equation}
\label{eq:interfoliation}
\lambda_{\fatA \min}\le\lambda_{{\fatA_\kappa} \min}\le \lambda_{{\fatA_\kappa} \max}\le\lambda_{\fatA \max}.
\end{equation}
{\bf Proof.} Reorder the rows and columns of $\fatA_\kappa$ such that the non-zero elements in the diagonal are in a symmetric block $\fatA_{\kappa\kappa}$ in the upper left corner of size $N_\kappa\times N_\kappa$.
The reordered and the original matrices have the same eigenvalues. The eigenvectors of $\fatA_\kappa$ with zero eigenvalues are non-zero in the upper $N_\kappa$ components and has one non-zero component in the lower part. 
The eigenvectors with non-zero eigenvalues have zeros in the lower $N-N_\kappa$ components and the upper part consists of the eigenvectors of $\fatA_{\kappa\kappa}$. The non-zero eigenvalues of $\fatA_\kappa$ coincide with the eigenvalues of $\fatA_{\kappa\kappa}$. 
Since $\fatA_{\kappa\kappa}$ is a principal submatrix on the diagonal of the symmetric $\fatA$, it follows from the eigenvalue interlacing property \cite[Ch. 8.1]{GVL} that the eigenvalues are ordered as in \eqref{eq:interfoliation}. 
$\blacksquare$
\vspace{2mm}

The time steps in the locally adaptive method in Algorithm \ref{alg:MRFE} are bounded by stability as follows
\begin{equation}
\label{eq:timestepbound}
  \Dtau_1\le \frac{2}{|\min_k \lambda_k(\fatA)|},\quad \Dtau_0=\frac{\Dtau_1}{m}.
\end{equation}
Then the stability condition \eqref{eq:forwardEstab2} for the short and long steps in \eqref{eq:multistab} is by the proposition and \eqref{eq:timestepbound}
\begin{equation}
\label{eq:timestepbound2}
\begin{array}{rl}
  |\Dtau_1\min_k\lambda_k(\fatA_1)|&= \displaystyle{\frac{2|\min_k\lambda_k(\fatA_1)|}{|\min_k \lambda_k(\fatA)|}\le 2},\\ |\Dtau_0\min_k\lambda_k(\fatA_0)|&=\displaystyle{\left|\frac{\Dtau_1}{m}\min_k\lambda_k(\fatA_0)\right|\le \frac{2|\min_k \lambda_k(\fatA_0)|}{m|\min_k \lambda_k(\fatA)|}<2}.
  \end{array}
\end{equation}
Both $\Dtau_0$ and $\Dtau_1$ in \eqref{eq:tau02def} satisfy a stability bound based on the eigenvalues as in \eqref{eq:forwardEstab2}.


\subsubsection{Estimate of work}\label{sec:estwork}

Assume that there are $p_0dN$ equations in $\calK_0$ and $p_1dN$ in $\calK_1$ with $p_\kappa\ge 0$ and $p_0+p_1=1$ in Algorithm \ref{alg:MRFE}. 
The computational work in each time step $\Dtau_\kappa$ is measured by the number of evaluations of the components of $\fatF$ ignoring the increased administration in the multirate method. 
With the same time step $\Dtau_0$ in all equations, the work is $m$ with $m$ evaluations of $\fatF$ and if the time step is $\Dtau_1$ then the work is 1. 

If longer time steps are taken in $\calK_1$, then the work for the whole system is $mp_0+p_1$. A part $p_0$ of $\fatF$ is evaluated $m$ times and a part $p_1$ once. The split system is preferred if 
\begin{equation}
\label{eq:workest}
   mp_0+p_1<m.
\end{equation}
If $m>1$, then
\[
   mp_0+p_1<mp_0+mp_1=m,
\]
and \eqref{eq:workest} is satisfied.

Let $|\eta_k|$ be ordered from large to small with growing $k$.
The distribution of $|\eta_k(t)|$ after proliferation at $t=0$ and $k=1$ is assumed to be
\begin{equation}
\label{eq:etaeq}
     |\eta_k(t)|=(\chi_0(t)-\chi_\infty)\exp(-(k-1)/k_e)+\chi_\infty.
\end{equation}
The estimate of $|\eta_k|$ for large $k$ is $\chi_\infty$, $\max_k|\eta_k|=\chi_0$ at $k=1$, $\chi_\infty\ll\chi_0$, and $k_e+1$ is the $k$ value where $\chi_0(t)-\chi_\infty$ is reduced by $1/e$.
According to \eqref{eq:chidefacc}
\[
     \chi_{a1}(t)=(\chi_0(t)-\chi_\infty)\exp(-(k_0-1)/k_e)+\chi_\infty=\frac{1}{m}\chi_{a0}(t),
\]
where $k_0$ is the number of equations in $\calK_0$. Since $\chi_\infty$ is small
\begin{equation}
\label{eq:k1}
   k_0\approx 1+k_e\ln m,
\end{equation}
independent of $t$. 
The number of equations in $\calK_1$ is $k_1\approx dN-1-k_e\ln m$. Thus, the quotients $p_0$ and $p_1$ in \eqref{eq:workest} are
\[
     p_0=\frac{1+k_e\ln m}{dN},\; p_1=1-\frac{1+k_e\ln m}{dN},
\]
and the work estimate for one time step $\Dtau_1$ after a proliferation is
\begin{equation}
\label{eq:workest2}
W_{\rm MRFE}=\displaystyle{mp_0+p_1}
            =\displaystyle{1+\frac{(m-1)(1+k_e\ln m)}{dN}}.
\end{equation}
The numerator of the last term in \eqref{eq:workest2} is independent of $N$. The number of evaluations of $\fatF$ approaches 1 per time step $\Dtau_1$ when $N$ is large, i.e. most equations are integrated with the largest time step. The alternative would be $m$ evaluations with the short time step $\Dtau_0$.
If $m=1$ then $\Dtau_0=\Dtau_1$ and one evaluation of $\fatF$ is required.

\subsection{Single rate time stepping with the backward Euler method}\label{sec:SRBE}

The ODE system \eqref{eq:poteq} is often stiff and it may be advantageous to solve it with an implicit method. A first order method is the backward Euler method \cite{HaNoWa}
to compute $\fatx^{n+1}$:
\begin{equation}
\label{eq:backwardE}
      \fatx^{n+1}=\fatx^n+\Dt \fatF(\fatx^{n+1}).
\end{equation}
The leading term in the local error is as in \eqref{eq:forwardEerr} with the opposite sign. 
The time step $\Dt$ for accuracy in this single rate method is chosen as in \eqref{eq:forwardEdt} for the forward Euler method. 
The backward method is stable for all systems with non-positive eigenvalues of the Jacobian and there is no time step restriction as in \eqref{eq:forwardEstab2}. 
This method has the same conservation properties as the forward 
Euler method in \eqref{eq:forwardEcons} applied to the gradient system \eqref{eq:poteq} and is frame invariant as in \eqref{eq:frameinv}.

The force in \eqref{eq:poteq} satisfies a one-sided Lipschitz condition. Use $\tfatA$ in \eqref{eq:forwardEdiff} to obtain
\begin{equation}
\label{eq:oneLip}
    (\fatu-\fatv)^T(\fatF(\fatu)-\fatF(\fatv))=(\fatu-\fatv)^T\tfatA(\fatu, \fatv)(\fatu-\fatv) \le c\|\fatu-\fatv\|_2^2,
\end{equation}
where $c$ is an upper bound on the eigenvalues of the symmetric part of $\tfatA(\fatu, \fatv)$ with $\fatu$ and $\fatv$ in a convex $\calB$
\[
      c=\max_{k, \fatu, \fatv} \lambda_k(\frac{1}{2}(\tfatA+\tfatA^T)),\; \fatu, \fatv\in \calB.
\]

The condition in \eqref{eq:oneLip} is invoked in the next theorem on nonlinear stability from \cite{StHu94} adapted for the backward Euler method. It is the discrete counterpart of Theorem 1 in Section~\ref{sec:gradsyst}.

\vspace{2mm}
\noindent
{\bf Theorem 3.}\hspace{2mm} The backward Euler method \eqref{eq:backwardE} for the gradient system in \eqref{eq:potential} with potential $V$ and ODE system \eqref{eq:poteq} satisfying the one-sided Lipschitz condition \eqref{eq:oneLip} is stable in the following sense
\begin{equation}
\label{eq:Vstab}
    V(\fatx^{n+1})-V(\fatx^n)\le\left(c-\frac{1}{\Dt}\right)\|\fatx^{n+1}-\fatx^n\|_2^2.
\end{equation}

\noindent
{\bf Proof.}\hspace{2mm} 
Use \eqref{eq:oneLip} and the mean value theorem as in \cite[Lemma 2.4.3]{HumphriesPhD} to prove
\begin{equation}
\label{eq:Vbound1}
    V(\fatu)-V(\fatv)\le\fatF(\fatu)^T(\fatv-\fatu)+ c\|\fatv-\fatu\|_2^2.
\end{equation}
Then by \cite[Res. 4.5]{StHu94}
\begin{equation}
\label{eq:Vbound2}
\begin{array}{rl}
    V(\fatx^{n+1})-V(\fatx^n)&\displaystyle{\le\left(-\frac{1}{\Dt}(\fatx^n-\fatx^{n+1})\right)^T(\fatx^n-\fatx^{n+1})+ c\|\fatx^{n+1}-\fatx^n\|_2^2}\\
       &\displaystyle{=\left(c-\frac{1}{\Dt}\right)\|\fatx^{n+1}-\fatx^n\|_2^2}.
\end{array}
\end{equation}
$\blacksquare$

\noindent{\bf Remark.}\hspace{2mm} 
We expect $c$ in \eqref{eq:Vstab} and \eqref{eq:Vbound2} to be small and positive in large time intervals and $\le 0$ close to steady state since it is close to $\max_k \lambda_k(\fatA(\fatx^n))$ (as it is in \eqref{eq:eigsmall}).

Using Theorem 3 and the lower bound on $V$ we can prove

\vspace{2mm}
\noindent
{\bf Corollary 1.}\hspace{2mm} Depending on $c$ in \eqref{eq:Vstab} let $\Dt$ satisfy
\begin{equation}
\label{eq:Dtbound1}
    c>0: \;\Dt\in(0, 1/c),\quad c\le 0:\; {\rm any}\; \Dt>0
\end{equation} 
Assume that the potential $V(\fatx^n),\, n=1,2,\ldots ,$ has a lower bound $V_{\min}$. The solution is determined by the backward Euler method in \eqref{eq:backwardE}. When $n\rightarrow\infty$, $V(\fatx^n)\rightarrow V^\infty$  and $\fatx^n\rightarrow \fatx^\infty$ where $V^\infty=V(\fatx^\infty)$ and $\fatx^\infty$ are constant and $V^\infty\ge V_{\min}$. If  $\fatx^n\in\calD$ and $\Dt$ is such that 
\begin{equation}
\label{eq:Dtbound2}
    \Dt\max_k\lambda_k(\fatA(\fatx^n))<1,
\end{equation} 
then $\fatx^{n+1}=\fatx^n=\fatx^\infty$.

\noindent
{\bf Proof.}\hspace{2mm} 
The potential $V(\fatx^n),\, n=1,2,\ldots,$ in Theorem 3 is non-increasing with $\Dt$ in \eqref{eq:Dtbound1}.
In addition, $V(\fatx^n)\rightarrow V^\infty\ge V_{\min}$ since $V$ has a lower bound and $\fatx^n\rightarrow \fatx^\infty$ since $c-1/\Dt<0$. 
The steady state $\fatx^\infty$ is such that $V(\fatx^\infty)=V^{\infty}$. If $\fatx^n\in\calD$ then $\fatF(\fatx^n)=\fatzero$ and one solution to \eqref{eq:backwardE} is $\fatx^{n+1}=\fatx^n=\fatx^\infty$. The Jacobian matrix $\fatI-\Dt\fatA(\fatx^n)$ is positive definite with the $\Dt$ in \eqref{eq:Dtbound2}. It follows from the implicit function theorem that this solution is unique and the claim is proved.
$\blacksquare$

\noindent{\bf Remark 1.}\hspace{2mm} 
Suppose that some distances $r_{ij}^\infty$ in $\fatx^\infty$ are such that $r_{ij}^\infty=s$ in \eqref{eq:gprop} and other distances satisfy $r_{ij}^\infty>r_A$. Then $g(r_{ij}^\infty)=0$ and $\fatF(\fatx^\infty)=\fatzero$. There are many such configurations of cells and 
the steady state $\fatx^\infty$ is not unique. It depends on the initial condition $\fatx^0$. 

\noindent{\bf Remark 2.}\hspace{2mm} 
The conditions in \eqref{eq:Dtbound1} and \eqref{eq:Dtbound2} are similar. If $\lambda_{\max}=\max_k\lambda_k(\fatA(\fatx^n))>0$ in \eqref{eq:Dtbound2} then $0<\Dt<1/\lambda_{\max}$ and if $\lambda_{\max}\le 0$ then $\Dt>0$. Moreover, $\lambda_{\max}\approx c$ in \eqref{eq:oneLip}. When the system approaches the steady state $\lambda_{\max}$ is small and positive. Close to steady state all eigenvalues are non-positive.

The system of nonlinear equations in \eqref{eq:backwardE} 
is solved numerically for $\fatx^{n+1}$ by a Newton-Krylov method \cite{BrHi86, KnollKeyes}. Define $\calF$ by
\begin{equation}
\label{eq:nonlinFdef}
      \calF(\fatx)=\fatx-\fatx^n-\Dt\fatF(\fatx)
\end{equation}
with the Jacobian matrix $\fatJ=\fatI-\Dt\fatA$. Solve 
\[
    \calF(\fatx^{n+1})=\fatzero
\]
for $\fatx^{n+1}$ by Newton iterations. Iterate as follows for  $j=0,1,\ldots$
\begin{enumerate}
\item  $\fatJ\Dfatx=-\calF(\fatx_{(j)}^{n+1})$,
\item  $\fatx_{(j+1)}^{n+1}=\fatx_{(j)}^{n+1}+\Dfatx$,
\end{enumerate}
and initialize with $\fatx^{n+1}_{(0)}=\fatx^n$ or $\fatx^{n+1}$ from \eqref{eq:forwardE}. 
The system of linear equations in the first step is solved for $\Dfatx$ by GMRES \cite{GMRES} using the explicit Jacobian $\fatA$. 
The outer Newton iterations are interrupted when 
\begin{equation}
    \Vert \Delta \textbf{x} \Vert_2 < \varepsilon_{\rm Newton} (\Vert \fatx_{(j+1)}^{n+1} \Vert_2 + 1), 
\end{equation}
where $\varepsilon_{\rm Newton}$ is chosen as $0.001\, \varepsilon$, or when the maximum number of outer iterations $n_{\rm Newton}$ is reached. Similarly, the GMRES iterations stop when 
\begin{equation}
\Vert \calF(\fatx_{(j)}^{n+1}) + \fatJ \Dfatx^{(k)}  \Vert_2 \leq \max(\veps_{\rm GMRES}\, \Vert \calF(\fatx_{(j)}^{n+1}) \Vert_2, \,\veps_{\rm GMRES}^{abs}),
\end{equation}
where we choose $\veps_{\rm GMRES}=\veps_{\rm GMRES}^{abs}= 0.001\, \veps$, or when the maximum number of iterations $n_{\rm GMRES}$ is reached \cite{scipy_sparse_linalg_gmres}. The algorithm is given in Algorithm \ref{alg:SRBE}.

\begin{algorithm}[tbp]
\linespread{1.15}\selectfont
\KwIn{Right-hand side \(\textbf{F}\) of Equation \eqref{eq:poteq},  start time \(t^0\), final time \(T\), initial coordinates \(\textbf{x}^0\), absolute accuracy \(\varepsilon\), Jacobian \(\textbf{A}\), maximum number of Newton iterations \(n_{\rm Newton}\), error threshold for Newton iterations \(\varepsilon_{\rm Newton}\)}
Initialize \(t = t^0\); \(\textbf{x} = \textbf{x}^0\)\;
\While{\(t<T\)}{
Evaluate force and Jacobian for current coordinates \(\hat{\textbf{F}} = \textbf{F}(\textbf{x})\); \(\hat{\textbf{A}} = \textbf{A}(\textbf{x})\)\;
Calculate time step size \(\Delta t= \sqrt{\tfrac{2\,\varepsilon}{\Vert \hat{\textbf{A}}\hat{\textbf{F}} \Vert_\infty }}\)\;
Initialize \(\textbf{x}_{next} =\textbf{x}\)\;
\For{\(j=0; j< n_{\rm Newton}; j++\)}{
\(\tilde{\textbf{F}} = \textbf{x}_{next} - \textbf{x} - \Delta t \hat{\textbf{F}} \)\;
\(\textbf{J} = I - \Delta t \hat{\textbf{A}}\)\;
Solve \(\textbf{J} \Delta \textbf{x} = - \tilde{\textbf{F}}\) for \(\Delta \textbf{x}\) using GMRES\;
Update \(\textbf{x}_{next} \leftarrow \textbf{x}_{next} + \Delta \textbf{x}\)\;
\If{\(\Vert \Delta \textbf{x}\Vert_2 < \varepsilon_{\rm Newton} (\Vert \textbf{x}_{next} \Vert_2 + 1)\)}{
\textbf{break}\;
}
Update force and Jacobian \(\hat{\textbf{F}} \leftarrow \textbf{F}(\textbf{x}_{next})\); \(\hat{\textbf{A}} \leftarrow \textbf{A}(\textbf{x}_{next})\)\;
}
Update \(\textbf{x} \leftarrow \textbf{x}_{next}\); \(t \leftarrow t + \Delta t\)\;
}
\SetAlgoRefName{IV}
\caption{Single rate backward Euler method (SRBE)}
\label{alg:SRBE}
\end{algorithm}

Assume that there is no proliferation and large time steps are possible to satisfy the accuracy requirements. The computational work per time step for the forward  Euler method SRFE in Algorithm \ref{alg:SRFE} and the backward Euler method SRBE in Algorithm \ref{alg:SRBE} are denoted by $\WfE$ and $\WbE$ and are measured by the number of evaluations of $\fatF$. 
The maximum time step for backward Euler is $\Dt_a$ given by accuracy and the maximum time step for forward Euler due to stability is $\Dt_s$. 
The time step for accuracy is the same for both forward and backward Euler.

Let $k_N$ be the number of steps in  the Newton iteration and let $k_G$ be the number of iterations in GMRES to solve \eqref{eq:backwardE}. 
The work to compute $\fatA$ and multiplication by $\fatA$ in GMRES is assumed to be about the same as one evaluation of $\fatF$. 
Thus, there are the equivalent to $k_{NG}=k_N(k_G+2)$ evaluations of $\fatF$ in the iterations in one time step of SRBE. Assuming $\fatF(\fatx^n)$ to be known, one evaluation is needed to estimate the local error in \eqref{eq:forwardEerrest} and one to compute $\fatF(\fatx^{n+1})$ in SRFE. Then
\[
         \WfE=2,\quad \WbE=k_{NG}+3=k_N(k_G+2)+3,
\]
and the total work in an interval $[0, T]$ with constant time steps is
\[
        \WfEtot=\frac{2T}{\Dt_s},\quad \WbEtot=\frac{(k_{NG}+3)T}{\Dt_a}.
\]
Then backward Euler is the best choice if $\WbEtot<\WfEtot$, i.e. if
\begin{equation}
\label{eq:workcmp}
      \frac{k_{N}(k_{G}+2)+3}{2}<\frac{\Dt_a}{\Dt_s}.
\end{equation}
As an example, take 3 Newton iterations and 3 GMRES iterations in each Newton step. Then backward Euler is more efficient if ${\Dt_a}/{\Dt_s}>9$. If $\Dt_a/\Dt_s<1$
then forward Euler is always the preferred method.

\section{Numerical results} \label{sec:numres}

In this section we compare the adaptive time stepping algorithms proposed in the previous section through numerical experiments with three cell configurations. We consider the following  adaptive algorithms: Algorithm \ref{alg:SRFE}, the single rate forward Euler method (SRFE); Algorithm \ref{alg:SRFES}, the single rate forward Euler method including the stability bound given by  \eqref{eq:forwardEstab2} (SRFES); Algorithm III, the two-level multirate forward Euler method (MRFE), and finally Algorithm IV, the single rate backward Euler method (SRBE). 
These algorithms are applied to the simulation of three different cell population configurations: (i) the relaxation between two daughter cells after division as a simple test case,
(ii) the relaxation of a spheroid in 3D where a single cell has been chosen to proliferate in the middle of the spheroid and (iii) a linearly growing tissue with multiple cell divisions of varying frequency as a biologically more realistic example. 


All experiments use the CBMOS package \cite{mathias2022cbmos}, extended with an implementation of the adaptive time stepping algorithms for the forward and backward Euler methods. The CBMOS software is a Python implementation of the center-based model, specifically designed for the numerical study of these models through the design of a flexible user interface exposing both forces and numerical solvers to the user. It is freely available under an MIT license on Github \cite{cbmos_github}. All Jupyter notebooks used to generate the figures in this section can be also be found in the repository.

The pairwise interaction force is the cubic force defined in \eqref{eq:cubic_force}.
The parameter values are chosen such that the time it takes for two daughter cells to relax to 99\% of the rest length \(s\) after proliferation (having been placed \(r_0=0.3\) cell diameters apart) corresponds to one hour (see \cite{MCBH20} for further details). Since the exact scaling of time is arbitrary, the time unit is this relaxation time \(\tau_{\text{relax}}\) and the length unit is the cell diameter or rest length $s$.  Table \ref{tab:parameters} lists the numerical parameters used throughout this section.

\begin{table}
\small
\renewcommand{\arraystretch}{1.25}
\centering
\begin{tabular}{c|l|r}
Parameter & Description & Value\\
\hline
\hline
\(s\) & rest length & 1.0 cell diameter \\
\(r_A\) & maximum interaction distance & 1.5 cell diameters\\
\(r_0\) & initial separation between daughter cells & 0.3 cell diameters\\
\(t^0\) & initial time& \(0\)\\
\(\mu\) & spring stiffness & 5.7\\ 
\(\epsilon\) & approximation parameter Jacobian-force product & 0.0001\\
\(\varepsilon\) & chosen absolute accuracy & 0.005 cell diameters \\
& & (if not specified differently)\\
\(m\) & ratio between levels for Algorithm III (MRFE) & 14 \\
& & (if not specified differently)\\
\(n_\text{Newton}\) & maximum number of Newton iterations & 5\\
\(\varepsilon_\text{Newton}\) & error threshold for Newton iterations &\(0.001\,\varepsilon\) \\ 
\(n_{\rm GMRES}\) & maximum number of GMRES iterations & 10\\
\(\varepsilon_{\rm GMRES}, \varepsilon_{\rm GMRES}^{abs}\) & error thresholds for GMRES iterations &\(0.001\,\varepsilon\) \\ 
\end{tabular}
\caption{Parameter values used in the numerical simulations throughout this section. In-simulation time is measured as multiples of the relaxation time between two daughter cells after proliferation $\tau_{\text{relax}}$. Length scales are measured in multiples of a cell diameter. }
\label{tab:parameters}
\end{table}

\subsection{Test case of the relaxation of two daughter cells after division}
We start by studying the time steps chosen by the single rate time stepping Algorithms I, II and IV for configuration (i)---two daughter cells relaxing after division---for different values of the chosen absolute accuracy \(\varepsilon\). For several reasons, this configuration represents an important test case. First of all, the CBM assumes that the forces acting on each cell can be expressed as a sum of pairwise interaction forces in \eqref{eq:poteqi}, making this a fundamental unit. Furthermore, as seen in Sections \ref{sec:forces_after_division} and \ref{sec:restr_step_size_proliferation}, the largest force magnitudes over the course of a simulation---requiring the smallest time steps to resolve accurately---occur right after cell division due to daughter cells being placed only a very short distance apart. Last but not least, it is possible to calculate the stability bound analytically for this simple case as \(\Delta t_s = {1}/{g'(s)}\). Note that for this simple test case all six equations evolve on the same time scale and we therefore do not consider the multirate Algorithm \ref{alg:MRFE} (MRFE).  

\begin{figure}
\begin{center}
\includegraphics[width=\linewidth]{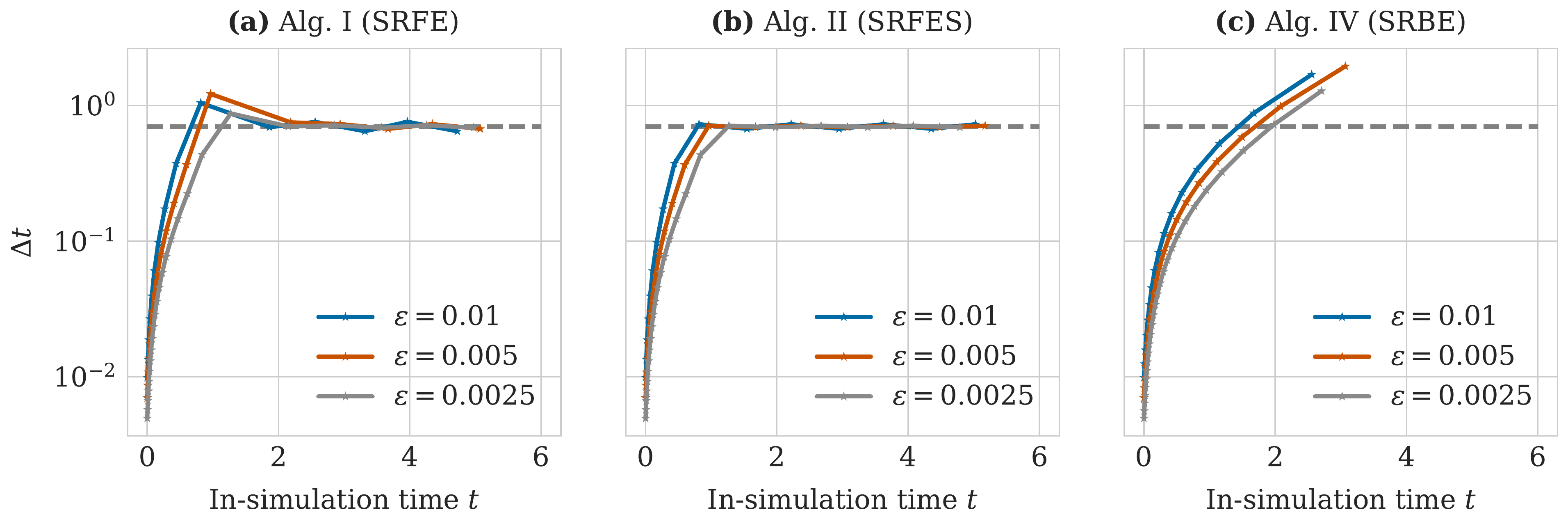}
\end{center}
\caption{Step size $\Delta t$ over time chosen adaptively for configuration (i), the relaxation of two daughter cells after proliferation, for different accuracy $\varepsilon$. The methods in the panels are: (a) Algorithm I; (b) Algorithm II; (c) Algorithm IV. The dashed horizontal line marks the analytical stability bound.
}
\label{fig:2cells_dt}
\end{figure}

Figure \ref{fig:2cells_dt} shows the time steps \(\Delta t\) used when simulating configuration (i) in $[0, 6]$ for $\veps= 0.01, 0.005, 0.0025$. 
In panel (a), Algorithm \ref{alg:SRFE} with the SRFE method calculates the trajectories of the midpoint coordinates. 
The time steps are determined by \eqref{eq:forwardEdt}. 
As the simulation progresses, the cells move apart and the magnitude of the pairwise force decreases, resulting in an increase of $\Dt$. Very soon $\Dt$ hits the stability bound and starts to oscillate around it. 
Note that the SRFE method does not explicitly take the stability limit into account. We observe that there is an overshoot in $\Dt$, yet the magnitude of the oscillations decreases in time. 
In panel (b), Algorithm \ref{alg:SRFES} chooses $\Dt$ according to \eqref{eq:forwardEdt} and explicitly calculates the stability bound \eqref{eq:forwardEstab2}. The large overshoot is removed but the oscillations for larger values of \(\varepsilon\) persist.
Lastly, in panel (c), the SRBE method from Section \ref{sec:SRBE} in Algorithm \ref{alg:SRBE}, is used. Since it is not limited by stability constraints, the step sizes increase in the time interval.

\begin{figure}
\begin{center}
\includegraphics[width=0.7\linewidth]{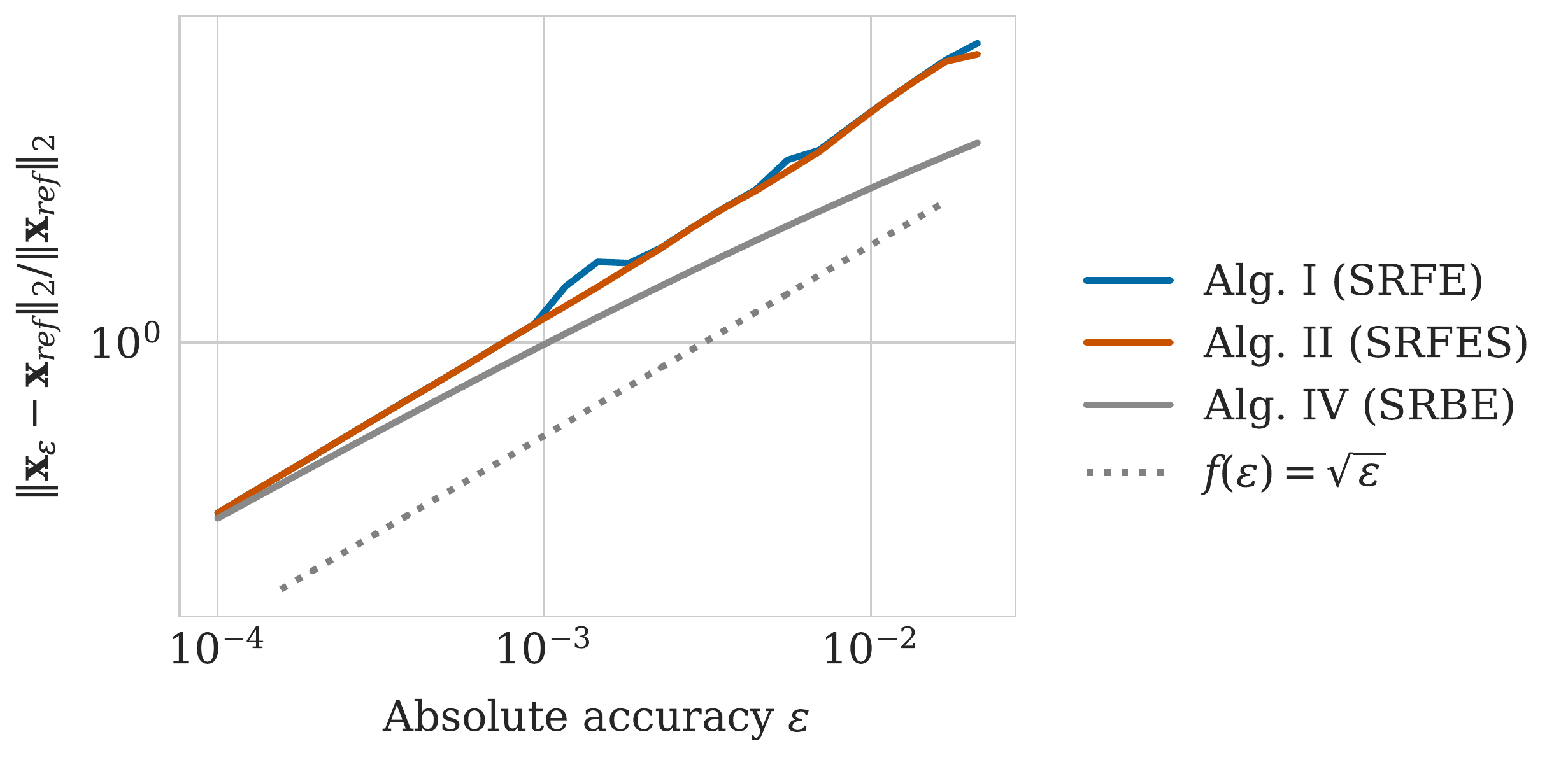}
\end{center}
\caption{Relative global error with respect to the reference solution as a function of $\varepsilon$ for configuration (i), two daughter cells relaxing after division, for the different time stepping algorithms. The performance is compared to the expected behavior in Theorem 2 (dotted line).} 
\label{fig:convergence_order}
\end{figure}

Next, we analyse the error with respect to a reference solution calculated with the forward Euler method and a fixed small time step \(\Delta t_{\rm ref} = 0.00005\). The cell trajectories calculated with the adaptive algorithms are interpolated down to the finer time resolution of the reference solution using a cubic interpolation scheme.
The relative error for $t\in[0, 3]$ is plotted as a function of \(\varepsilon\) in Figure \ref{fig:convergence_order}. The error in space is measured  in \(\Vert \cdot \Vert_{\infty}\) 
and then the \(\Vert \cdot\Vert_2\) norm is applied over the resulting time series.
For the SRFE and SRFES methods, our results confirm that the error decays as $\sqrt{\veps}$ as expected from Theorem 2 with $p=1$. 
For the SRBE method  in Algorithm IV we obtain a numerical order of convergence of 0.42 using least squares regression. 

\subsection{Single proliferation event within a larger spheroid}

While looking at two cells in isolation represents a convenient numerical test case, in more realistic settings cells divide with neighbors around them which contribute to the overall forces experienced. It is therefore of interest to consider cell proliferation within larger cell populations.
To this end, we let a single division event take place within a 3D spheroid of 216 cells (six cells in each dimension), where the cell midpoints have been arranged on a hexagonal close packed lattice configuration (configuration (ii)). The distance between neighboring cells is chosen as exactly one rest length $s$, so that no forces are active between them. Before the start of the simulation the middle cell in the spheroid is deleted and two daughter cells are placed at \(r_0=0.3\) apart with the cell division in a random direction, such that the midpoint between them is situated at the former position of the mother cell. The seed of the random number generated is fixed across all numerical experiments for reproducibility, resulting in the same cell division direction being drawn for all of them. 

For this test case the locally adaptive multirate method, Algorithm \ref{alg:MRFE}, is of special interest. We hence consider it in addition to the single rate methods Algorithms I, II and IV. 

\subsubsection{Study of time step sizes}
\begin{figure}
\begin{center}
\includegraphics[width=0.75\linewidth]{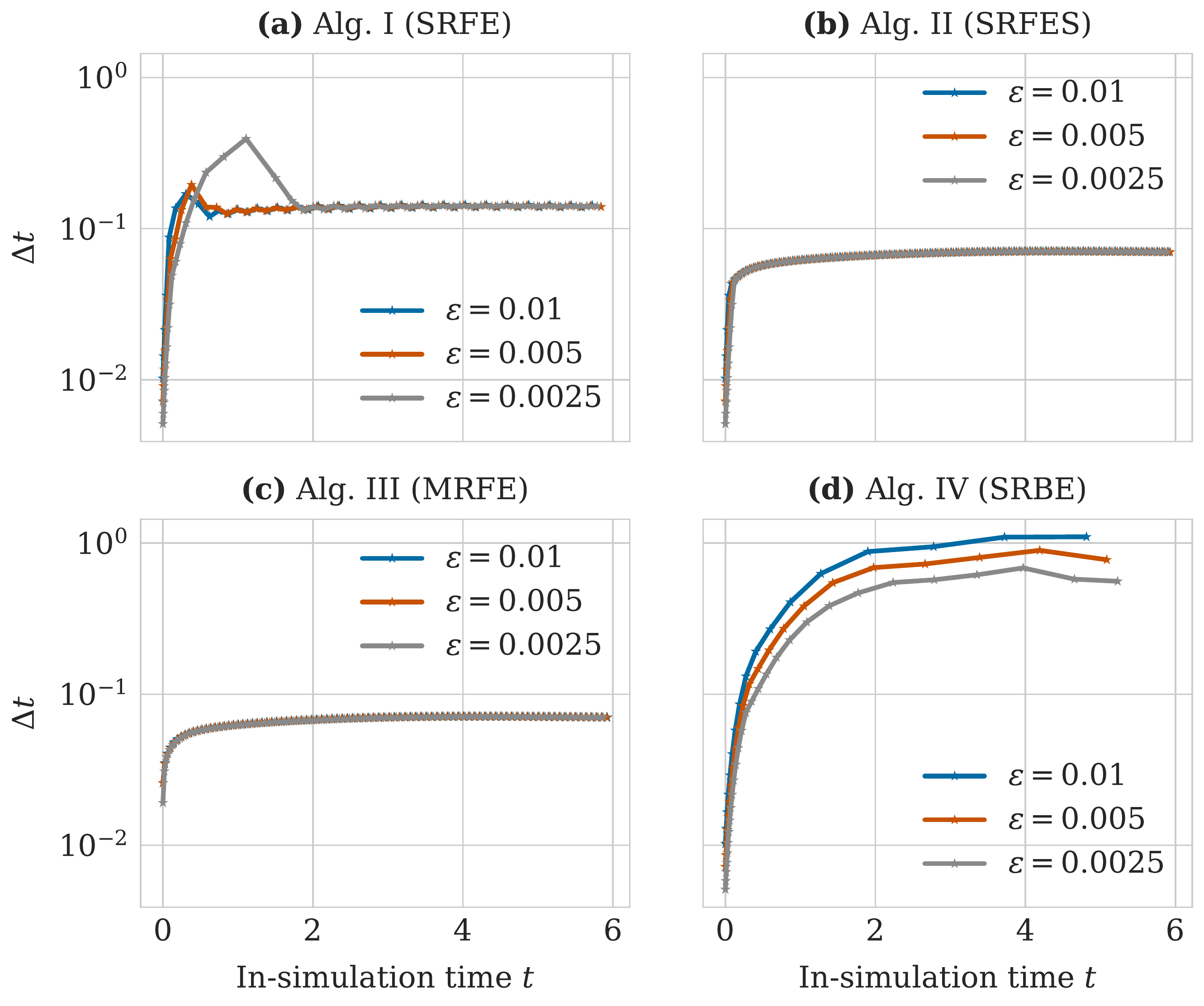}
\end{center}
\caption{Step size $\Delta t$ over time chosen adaptively for a single proliferation event taking place within a spheroid of cells for different absolute accuracy values $\varepsilon$. The panels differ in the time stepping algorithm chosen: (a) single rate forward Euler method; (b) single rate forward Euler method with stability check; (c) multirate forward Euler method, here the largest time step used \(\Delta t_\text{max}=\Dtau_1\) is plotted; (d)  single rate backward Euler method. }
\label{fig:prolif_within_spheroid_dt}
\end{figure}

We again start by studying the time steps chosen by the different adaptive algorithms for different values of \(\varepsilon\). 
The time steps calculated for the simulation of configuration (ii)
are shown in Figure \ref{fig:prolif_within_spheroid_dt} with $t\in[0, 6]$ for a randomly chosen seed that we fix across algorithms and error tolerances. 
The seed 
affects the exact positions of the daughter cells and resulting pairwise forces with surrounding neighbors. Nevertheless, the variation in force magnitude depending on the exact cell division direction is small and the time step sizes can be expected to be representative. 

In panel (a) we observe that similarly to configuration (i), the time steps in Algorithm I are initially restricted by accuracy and then quickly increase until being restricted by stability before $t=2$. Again, we observe for small $t$ an overshoot that increases in magnitude for smaller values of \(\varepsilon\) and oscillations for larger $t$.  The difference between the initial time step  size \(\Delta t_0\) and the time step size at steady state \(\Delta t_s\), here given by the stability bound, is about one order of magnitude for \(\varepsilon=0.01\) and \(\varepsilon=0.005\) and roughly 1.5 orders of magnitude for  \(\varepsilon=0.0025\). 

The time steps calculated by Algorithm II in panel (b) use the same initial time steps \(\Delta t_0\), but do not exhibit overshoot as it takes the stability bound into account. Furthermore, the time step at steady state \(\Delta t_s\)  is smaller. This is due to the fact that the eigenvalues of the Jacobian needed for stability are underestimated using Gershgorin's theorem \eqref{eq:Gersh}. It is not analytically possible to calculate the stability limit for this case and that calculating the eigenvalues fully is too computationally costly. 

Algorithm III with MRFE also uses Gershgorin's theorem to estimate the eigenvalues and calculate the stability bound. In Figure \ref{fig:prolif_within_spheroid_dt} (c), $m=14$ and the largest time step  $\Dtau_1$ is plotted. Consequently, the initial time step is larger than the one used by the Algorithms I and II in panels (a) and (b). A small number of equations are initially solved using multiple small time steps $\Dtau_0=\Dtau_1/m$. It is necessary to impose an explicit stability bound in the MRFE method. Otherwise, the time steps will oscillate wildly.

The time step sizes chosen by Algorithm IV are found in Figure \ref{fig:prolif_within_spheroid_dt}(d). At the beginning of the simulation, the time steps are the same as those determined by Algorithms I and II in panels (a) and (b) and restricted by accuracy. 
They stabilize once the system has reached a steady state after relaxation of all forces at about $t=2$ and are bounded by accuracy. They are
an order of magnitude larger than both the correct stability limit found by the adaptivity in Algorithm I and the Gershgorin estimate used by Algorithms II and III. 

\begin{figure}
\begin{center}
\includegraphics[width=\linewidth]{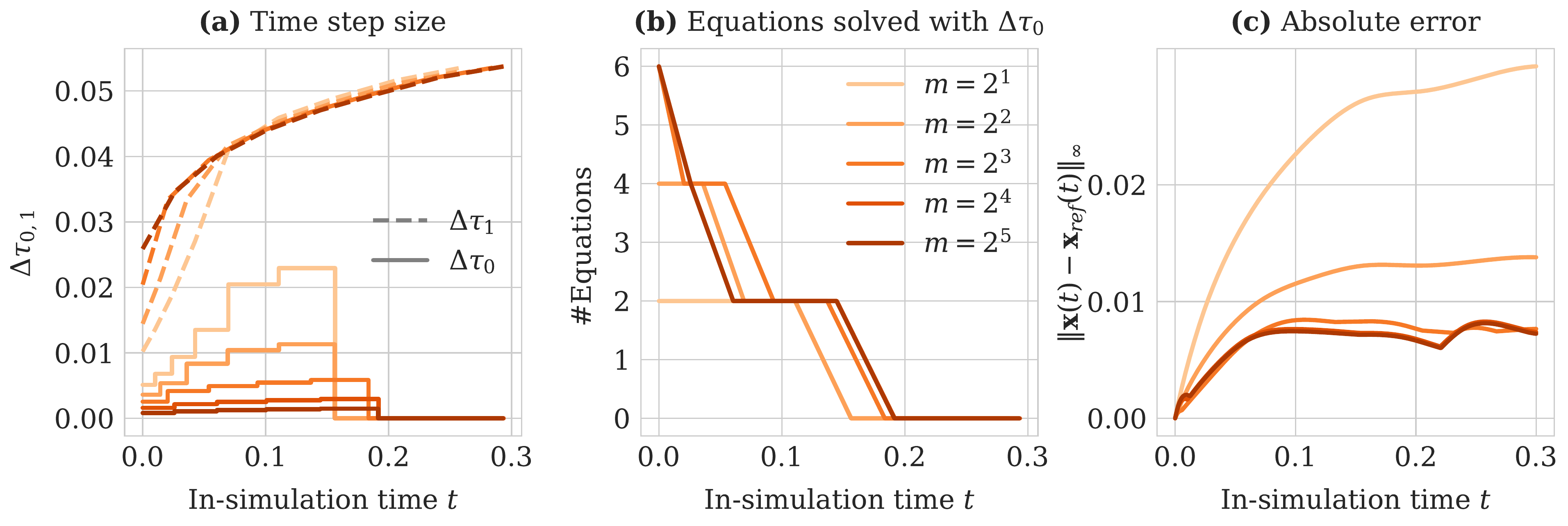}
\end{center}
\caption{The dependence of time steps and number of equations per level used by MRFE on the ratio \(m\) between levels when applied to configuration (ii) with \(\varepsilon=0.005\). (a) Time steps for both levels over time. Solid lines refer to \(\Delta \tau_0\), dashed lines to \(\Delta \tau_1\). The colors are explained by the legend in panel (b). (b) Number of equations solved with multiple small time steps \(\Delta \tau_0\) over time. (c) The error for different $m$.}
\label{fig:comp_dt_MRFE}
\end{figure}

The dependence of the behavior of the multirate Algorithm III on the ratio between levels $m$ in \eqref{eq:taudef} is displayed in Figure \ref{fig:comp_dt_MRFE}. In panel (a), we take the time steps $\Dtau_0$ and $\Dtau_1$ for different ratios \(m\), \(\varepsilon=0.005\), and time $t$. After a full time step \( \Dtau_1\), the time step on both levels is increased. The number of equations solved with \(\Dtau_0\) are plotted in panel (b). 
We observe that Algorithm III initially chooses two levels with the equations of the two daughter cells on the fine level.
After only a few time steps with \(\Dtau_1\), there is a single level independent of \(m\). 
The error in \(\Vert\cdot\Vert_\infty\) with respect to a fixed time step reference solution calculated with \(\Delta t_{\rm ref} = 0.0005\) is shown in Figure \ref{fig:comp_dt_MRFE} (c). The error decreases for ratios \(m\ge 4\). 

\begin{figure}
\begin{center}
\includegraphics[width=0.45\linewidth]{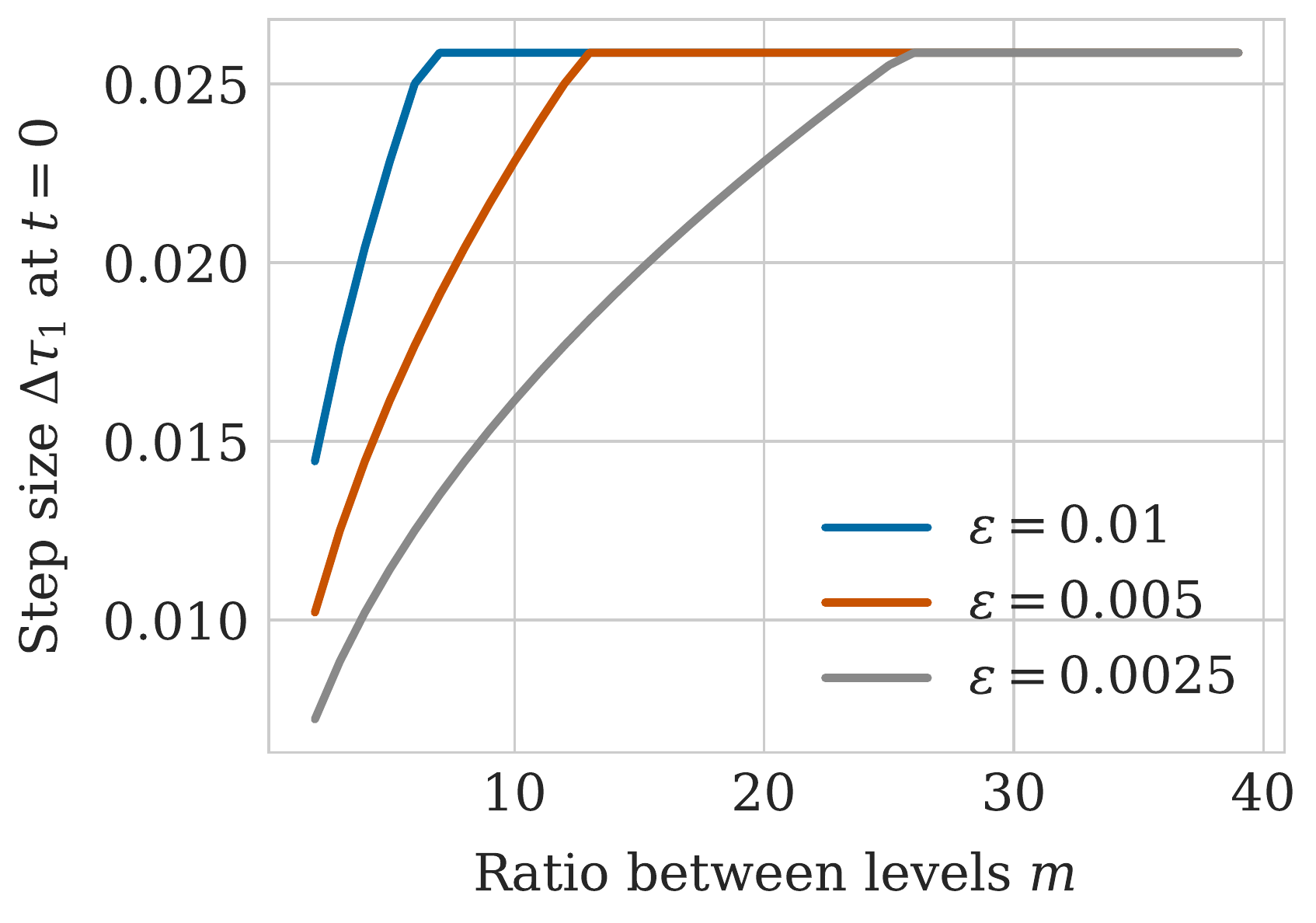}
\end{center}
\caption{Time step size \(\Dtau_1\) for larger level in MRFE at \(t=0\) for different absolute accuracy values \(\veps\) as a function of the ratio \(m\).}
\label{fig:optimal_m}
\end{figure}

Figure \ref{fig:optimal_m} plots the size of the initial time step on the coarser level \(\Dtau_1\) as a function of the ratio \(m\) for different absolute accuracy values \(\veps\). We observe that there is a limit on the size of  \(\Dtau_1\) at \(t=0\) independently of \(m\). This limit consists of the stability limit \(\Delta t_s\) and is reached for different \(m\) values depending on the chosen accuracy. We can conclude that the optimal \(m\) for a given accuracy is the smallest \(m\) such that the equations on the coarser level are solved with the stability time step \(\Delta t_s\), as increasing \(m\) further increases the number of partial updates necessary, while not allowing for larger step sizes. 
In Figure \ref{fig:prolif_within_spheroid_dt} as well as the following numerical experiments,  let the ratio between the two levels be \(m=14\) which is optimal for \(\veps=0.005\). 

\begin{figure}[tbp]
\begin{center}
\includegraphics[width=\linewidth]{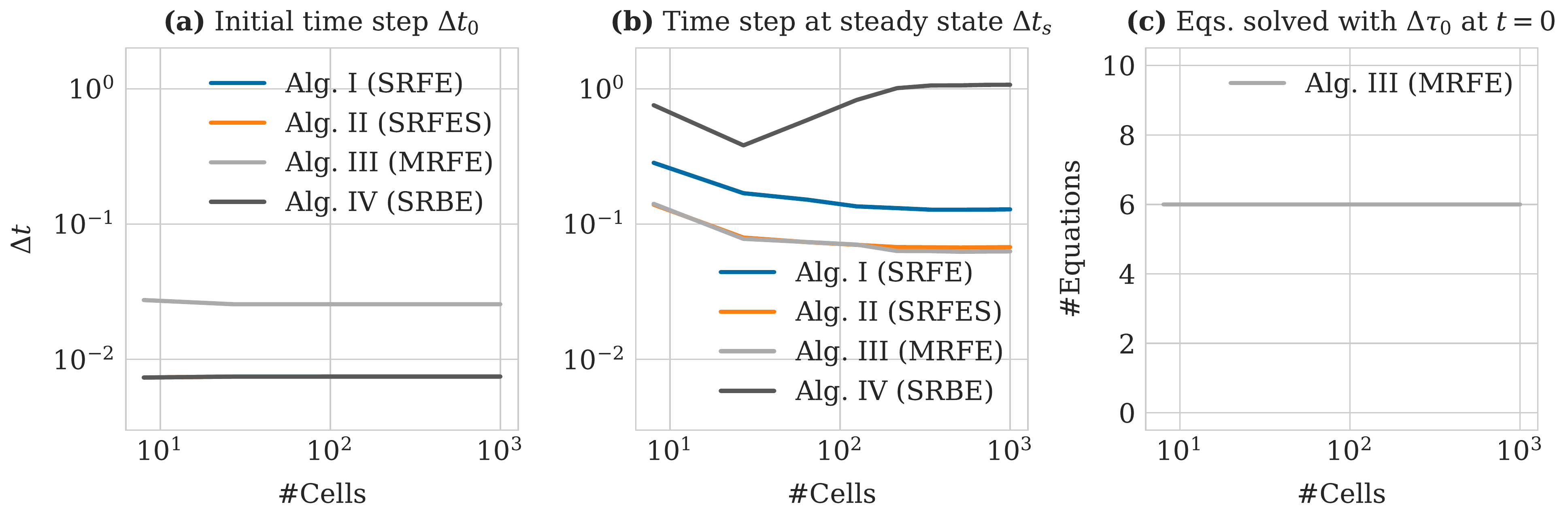}
\end{center}
\caption{Dependence of time steps and number of equations initially solved on the number of cells $N$ in a spheroid. (a) Size of the initial time step $\Delta t_0$ after proliferation as a function of $N$ for four different algorithms. The curves for Algorithms I, II, and IV overlap as they take the exact same initial time step. (b) Size of time step at steady state $\Dt_s$ as a function of the number of cells $N$. The curves for Algorithms II and III overlap, since they estimate the same stability bound. (c) Number of equations initially solved using \(\Delta t_\text{min}=\Dtau_0\) in the MRFE method as a function of the number of cells.}
\label{fig:dependence_n_cells}
\end{figure}

In the next experiment, we are interested in studying the dependence on the number of cells of the initial step size \(\Delta t_0\) and the step size at steady state \(\Delta t_s\) used by Algorithms I-IV. To this end, we consider a range of cell spheroids of different sizes and choose the accuracy parameter \(\varepsilon = 0.005\). The number of cells $N$ varies between 8+1 to 1000+1 cells (the additional cell is due to the second daughter cell) in Figure \ref{fig:dependence_n_cells}. 
For each spheroid size, we averaged the chosen $\Dt$ over five different initial cell division directions to remove any effect of the initial placement. 

As we can see from panel (a), the initial time step is independent of the number of cells in the spheroid as expected in Section \ref{sec:restr_step_size_proliferation}. 
The daughter cells only experience forces after proliferation from their direct neighbors and hence it is irrelevant if the spheroid contains 8 or 1000 cells. The minimum time steps are similar when $N=2$ in Figure \ref{fig:2cells_dt}. For Algorithm III, we plot the largest time step, which is why the initial time step is larger than for Algorithms I and II. 
From (b) the conclusion is that the time step at steady state, i.e.\ the stability bound $\Dt_s$ for Algorithms I-III, and the steady state time step $\Dt_a$ bounded by accuracy for Algorithm IV, decreases slightly for larger spheroid sizes. Since $\Dt_a$ is the time step also for the forward Euler methods to satisfy the accuracy requirement, the quotient $\Dt_a/\Dt_s$ in \eqref{eq:workcmp} is about ten for Algorithms I-III.
The stability bound and the minimum of the eigenvalues of the Jacobian are almost independent of $N$ as concluded in the end of Section \ref{sec:explstab}. 
Finally in panel (c), the number of equations initially solved with multiple smaller time steps in Algorithm III is also independent of the size of the spheroid
because the strongest forces in the cell system are between the two new daughter cells with six degrees of freedom, see Section \ref{sec:forces_after_division}.

\subsubsection{Convergence results}

\begin{figure}
\begin{center}
\includegraphics[width=0.8\linewidth]{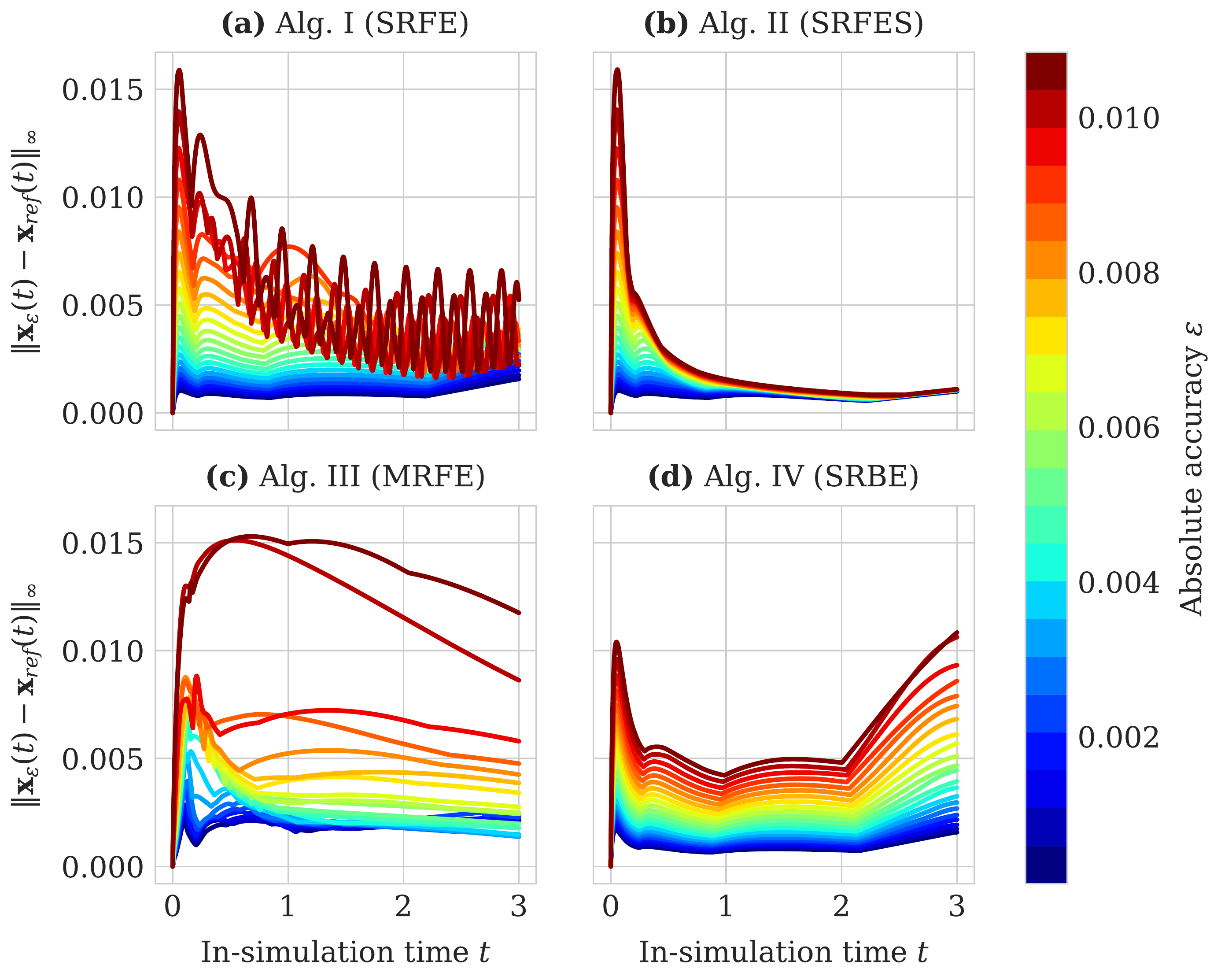}
\end{center}
\caption{Absolute difference in \(\Vert\cdot\Vert_\infty\) to a reference solution $\fatx_{\rm ref}$ for different values of $\varepsilon$ for configuration (ii) with
(a) Algorithm I; (b) Algorithm II; (c) Algorithm III; (d) Algorithm IV. 
}
\label{fig:convergence_error_spheroid}
\end{figure}

The cell trajectories calculated by the adaptive Algorithms I-IV are compared to a reference solution $\fatx_{\rm ref}$ in Figure \ref{fig:convergence_error_spheroid}. We measure their absolute difference (approximately the global error $\fatE$ in Theorem 2) in \(\Vert\cdot\Vert_\infty\) for $t\in[0, 3]$ and a large range of \(\varepsilon\). The reference solution $\fatx_{\rm ref}$ uses the fixed time stepping forward Euler method with \(\Delta t_\text{ref} = 0.0005\). The coarser solutions are interpolated down to the finer time grid using a cubic interpolation scheme. 

We observe oscillations for Algorithm I (Figure \ref{fig:convergence_error_spheroid} (a)) for large values of \(\varepsilon\) as it does not explicitly check the stability bound. All other algorithms have smooth error curves for all \(\varepsilon\) values. Algorithm II (Figure \ref{fig:convergence_error_spheroid} (b)) displays errors quickly decreasing over time even for large  \(\varepsilon\) values and as such has the smallest errors overall. In comparison, Algorithm III (Figure \ref{fig:convergence_error_spheroid} (c)) has larger errors for all \(\varepsilon\) values, but still on the order of \(\varepsilon\). 
Algorithm IV in Figure \ref{fig:convergence_error_spheroid} (d) displays errors comparable to Algorithm I albeit without the oscillations for large \(\varepsilon\) values. Also, the errors do not decay over time as they do for Algorithm II, but stay constant or even increase for large \(\varepsilon\) values. For small $\veps$, the algorithms display very similar errors. 


\begin{figure}
\begin{center}
\includegraphics[width=\linewidth]{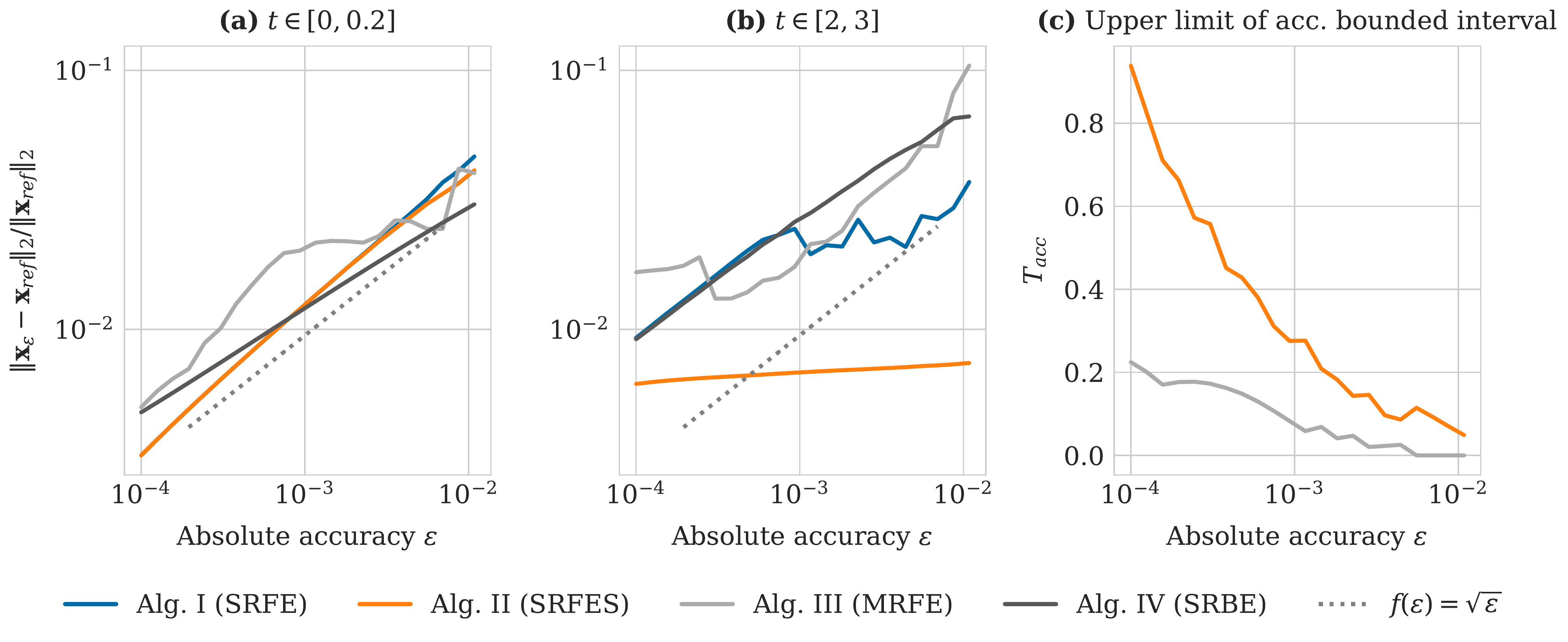}
\end{center}
\caption{(a), (b) Relative error with respect to the reference solution $\fatx_{\rm ref}$ as a function of $\varepsilon$ for configuration (ii)
with Algorithms I-IV. 
The time interval is in (a) \([0, 0.3 ]\)
and in (b) \([2, 3]\).   
(c) Upper limit of the accuracy bounded interval \(\Tacc\) as a function \(\varepsilon\) for Algorithms II and III. 
}
\label{fig:convergence_order_spheroid}
\end{figure}

The relative error in the solution with the four algorithms is found in Figure \ref{fig:convergence_order_spheroid}. In panel (a), the time interval is $[0, 0.3]$. The accuracy determines $\Dt$ and the error is proportional to $\sqrt{\veps}$ in the explicit methods in accordance with Theorem 2. 
The stability bounds $\Dt$ in $[2, 3]$ with the stability check in Algorithms II and III in panel (b) with an error almost independent of $\veps$. When $\veps<10^{-3}$ in Algorithm I and for all $\veps$ in the implicit Algorithm IV, accuracy puts a limit on $\Dt$ and the error decays as $\sqrt{\veps}$. This is in agreement with the estimates in \eqref{eq:Aprolifrestr} after proliferation. The accuracy determines $\Dt$ in $[0, \Tacc]$ and the stability when $t>\Tacc$ for Algorithms II and III in panel (c). The stricter the accuracy requirements are, the longer the interval is for both methods.  

\subsubsection{Comparison of computational cost}
\label{subsubsec:comp_cost}

In this section, we evaluate the benefit of using an adaptive time stepping algorithm for solution of the cell trajectories in configuration (ii) in terms of computational cost, both in the number of force function $\fatF$ and Jacobian $\fatA$ evaluations and in average wall clock time elapsed to simulate until a certain in-simulation time \(T\). 
As a baseline we use the forward Euler method with a fixed time step chosen as the initial time step $\Dt_0$ in Algorithms I and II. 
This is the step size that is necessary to resolve the cell trajectories everywhere to the desired absolute accuracy of \(\varepsilon=0.005\) in $[0, T]$ and is in general not known \emph{a priori}. The fixed step method computes $\fatF$ once in each time step.


A larger spheroid is created with \(13^3+1=2198\) cells with the same seed for all algorithms to randomize the direction of the initial cell division. The wall clock times are averaged over 20 repetitions. The benchmarks were run on Rackham, a high performance cluster provided by the Multidisciplinary Center for Advanced Computational Science (UPPMAX). The node we used consisted of two 10-core Intel Xeon E5-2630 v4 processors at 2.2~GHz, 128~GB of memory.

\begin{figure}
\begin{center}
\includegraphics[width=\linewidth]{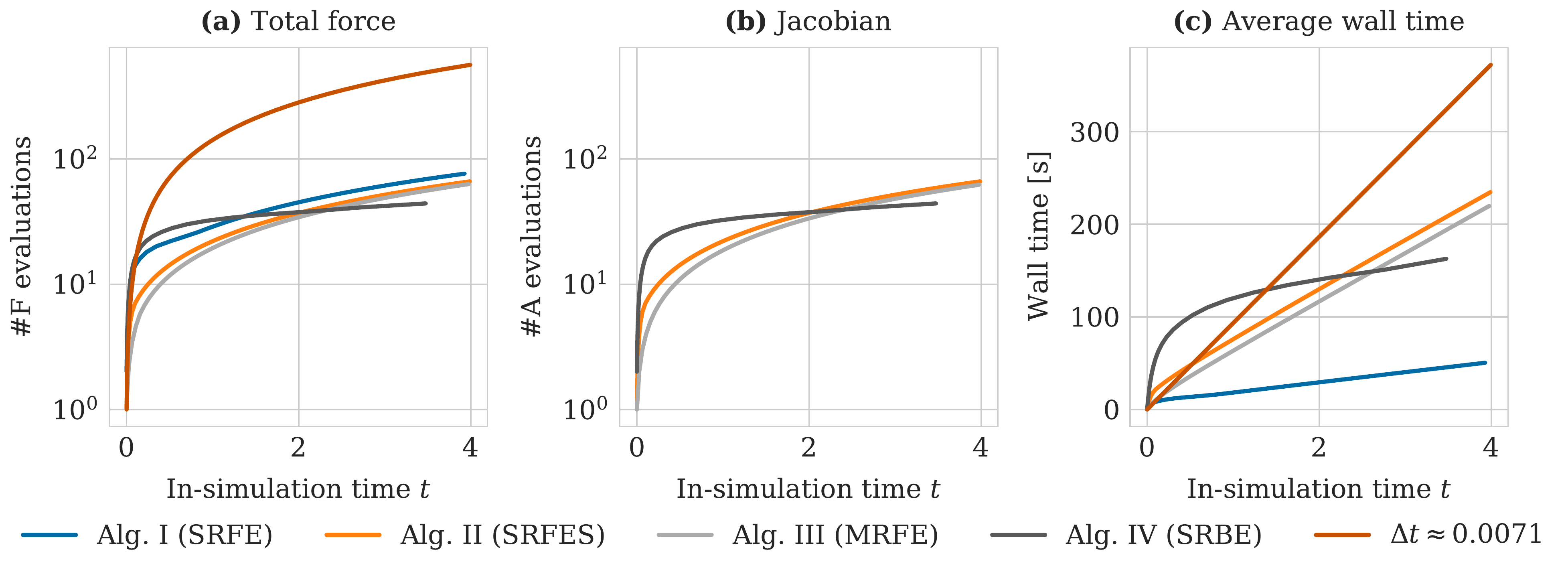}
\end{center}
\caption{Comparison of computational costs for Algorithms I-IV with a forward Euler method with a fixed time step \(\Delta t_{\rm fixed} \approx 0.0071\) for simulation of configuration (ii), a single initial proliferation event within a larger spheroid. (a, b) Number of (a) total force \(\textbf{F}\) and (b) Jacobian \(\textbf{A}\) evaluations as a function of in-simulation time \(t\).  (c) Wall time averaged over 20 repetitions and a fixed random seed as a function of in-simulation time \(t\). }
\label{fig:benchmark_single_prolif_spheroid}
\end{figure}

Figure \ref{fig:benchmark_single_prolif_spheroid} shows the number of total force evaluations in panel (a) as function of the in-simulation time \(t\). Partial updates for Algorithm III (MRFE) are counted as the number of perturbed equations affected by the partial update divided by the total number of equations as in Section \ref{sec:estwork}. The method with fixed time steps requires significantly more evaluations of the forces $\fatF$ than the adaptive methods in panel (a). Out of those, Algorithms II and III require less evaluations than Algorithm I as they do not approximate the Jacobian-force product using an additional force evaluation as in \eqref{eq:forwardEerrest}. Instead they evaluate the Jacobian as seen in panel (b). 
Note that Algorithm I (SRFE) and the fixed time stepping forward Euler method do not require evaluations of the Jacobian \(\textbf{A}\). Algorithm IV initially requires both more force and Jacobian evaluations compared to Algorithms I-III. However, this ratio shifts once the step size has been limited by stability for those methods for a sufficient number of steps and in-simulation time \(t\). 

Panel (c) of Figure  \ref{fig:benchmark_single_prolif_spheroid} shows the  wall clock time again as a function of the in-simulation time \(t\). Here, Algorithm I, the SRFE method with an adaptive time step without a stability check, is the fastest method. The reason for the efficiency of the SRFE method appears to be that $\Dt$ is chosen adaptively and that there is no calculation of the smallest eigenvalue of the Jacobian to impose a bound on $\Dt$. Moreover, there is less administration in SRFE in every time step. The disadvantage of SRFE is the larger and more oscillatory errors in Figure \ref{fig:convergence_error_spheroid}(a).
Algorithms II and III are both more efficient than the fixed time stepping for a longer simulation interval, with the multirate Algorithm III being marginally better than the single rate Algorithm II.
Algorithm IV is more competitive in the later half of the interval, when the time step of the explicit methods is bounded by stability. 

Using the data in Figures  \ref{fig:benchmark_single_prolif_spheroid} (c) and \ref{fig:dependence_n_cells} (b), the discussion in the end of Section \ref{sec:SRBE}, and \eqref{eq:workcmp}, the observed relation between the work in one time step with Algorithms IV (SRBE) and I (SRFE) after the transient phase with $t>1$ is
\[
   \frac{\WbE}{\WfE}=\frac{\Dt_a}{\Dt_s}\frac{\WbEtot}{\WfEtot}\approx 10\cdot 1=10.
\]
Algorithm IV is considerably more demanding than Algorithm I in terms of work in one time step but the total work for Algorithm IV when $t>2$ is comparable in time and lower in the number of $\fatF$ and $\fatA$ evaluations.

The work after proliferation measured in wall clock time for the methods in Figure \ref{fig:benchmark_single_prolif_spheroid}(c) is approximately
\begin{equation}\label{eq:workapp}
    W(t)=W_0(1-\exp(-t/\tau))+W_1 t, 
\end{equation}
where $W_0, W_1,$ and the time scale for the transient $\tau$ are constant. 
First suppose that the proliferations occur regularly at $t=(i-1)\DT,\, i=1,2,\ldots \nu,$ in an interval $[0, T]$ with $T=\nu \DT$. 
Then the total work is by \eqref{eq:workapp}
\begin{equation}\label{eq:workdet}
    \sum_{i=1}^{\nu}W(\DT)=\sum_{i=1}^{\nu}W_0(1-\exp(-\DT/\tau))+W_1\DT=\nu W_0(1-\exp(-T/(\nu\tau))+W_1 T. 
\end{equation}
Now suppose that the number of proliferations in $[0, T]$ is Poisson distributed. The interval between the proliferations is $\DTs$ and it is exponentially distributed with rate parameter $\lambda$. The probability density function for $\DTs$ is $\lambda\exp(-\lambda\DTs)$. The expected $\DTs$ is $1/\lambda$ and let the rate be $\lambda=\nu/T$. Then the expected total work is
\begin{equation}\label{eq:workstoch}
    \Expc[\sum_{i=1}^{\nu}W(\DTs)]=\sum_{i=1}^{\nu}W_0\Expc[(1-\exp(-\DTs/\tau))]+W_1\Expc[\DTs]=\frac{\nu W_0}{1+\nu\tau/T}+W_1 T. 
\end{equation}
Compare the total work in \eqref{eq:workdet} and \eqref{eq:workstoch}. Since $1/(1+\nu\tau/T)$ is an approximation of $1-\exp(-T/(\nu\tau))$, the expected total work with a random proliferation is well approximated by the total work with a deterministic proliferation.


\subsection{Linearly growing spheroid}

As a final numerical experiment we consider the case of linear tissue growth where the cells divide at a (globally) fixed frequency. We start with a spheroid of \(N_0 = 2197\) and let it increase its number of cells by \(n=10\) to \(N_T = 2207\). For simplicity, we deterministically choose the time between two cell division events at the population level to be \(\Delta t_{\rm div}\), and generate the set of cell division times \(\{i\, \Delta t_{\rm div} | i= 1,...,n\}\) accordingly in advance. During the simulation, we choose a random cell at each cell event time and let it divide into two daughter cells. The final simulation time is given by \(T = n\, \Delta t_{\rm div}\). We measure the wall time it takes to solve this using the adaptive Algorithms \ref{alg:SRFE}-\ref{alg:SRBE} and an absolute accuracy \(\veps=0.005\). As a baseline, we again use a fixed time stepping algorithm with \(\Delta t_{\rm fixed}= 0.0078\) given by \eqref{eq:Aprolifrestr}. 
The wall times of Algorithms I-IV are measured as multiples of the baseline wall time as follows. For each data point the absolute wall times are averaged over four repetitions with the same random seed and then divided by the baseline wall times. 
All absolute wall times are found in Table \ref{tab:linear_growth}. 
As before we run the simulations on the Rackham high performance cluster (see Section \ref{subsubsec:comp_cost}).

\begin{figure}
\begin{center}
\includegraphics[width=0.45\linewidth]{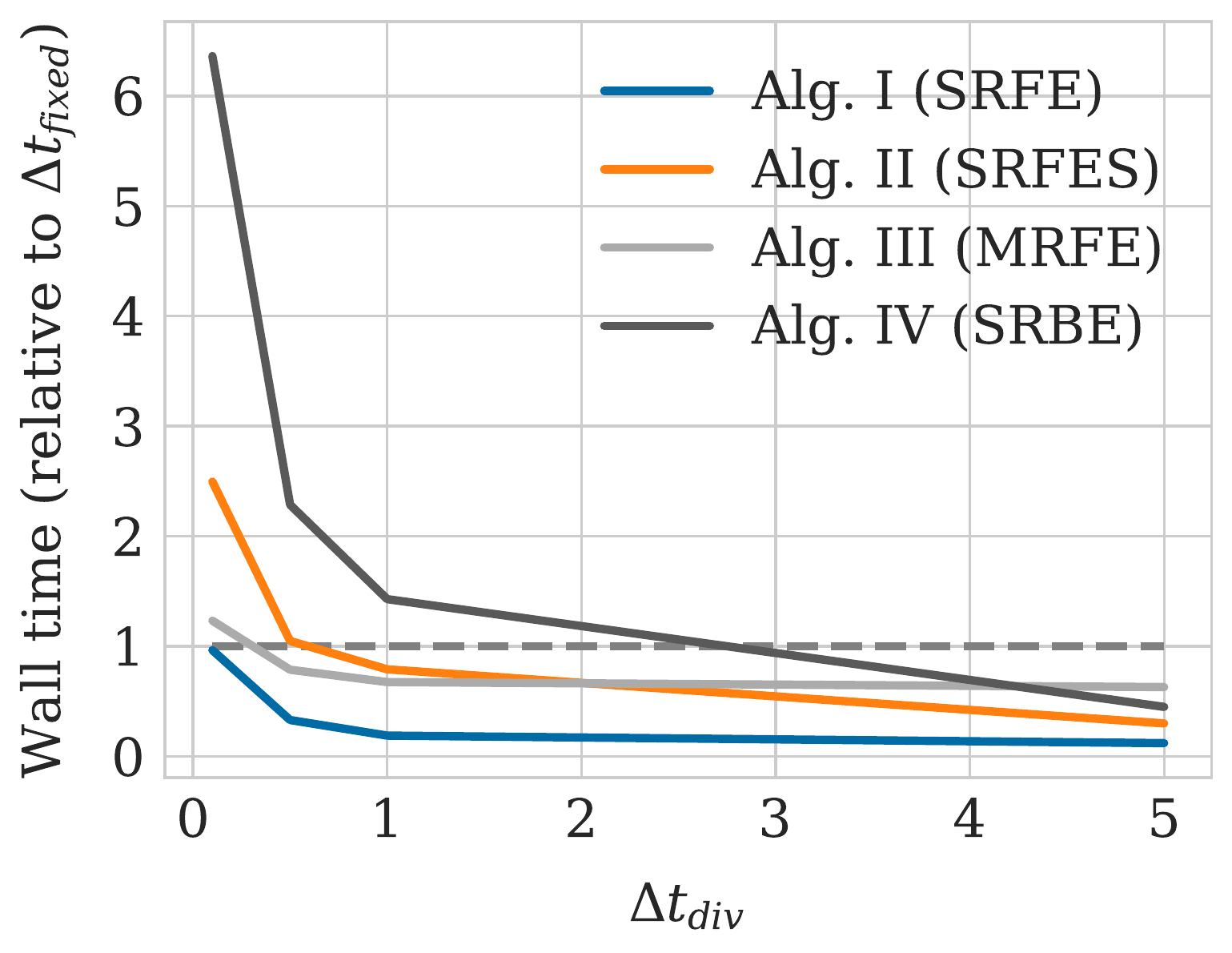}
\end{center}
\caption{Relative wall times for linear tissue growth as a function of time between cell division events \(\Delta t_{\rm div}\) on the population level for the Algorithms I-IV for configuration (iii). 
Wall times are calculated relative to the wall time for the fixed time stepping and averaged over four repetitions. The dashed line visualizes wall times equal to the fixed time stepping.
The end simulation times \(T\) and the absolute wall times for all algorithms can be found in Table \ref{tab:linear_growth}.}
\label{fig:linear_tissue_growth}
\end{figure}

Figure \ref{fig:linear_tissue_growth} displays the relative wall time of Algorithms I-IV as a function of the time between consecutive cell division events \(\Delta t_{\rm div}\). Algorithm I (SRFE) is beneficial compared to fixed time stepping for all values \(\Delta t_{\rm div}\) considered. Unless cell proliferation is very frequent, i.e.\ happening at a global scale more often than 5 times the relaxation time \(\tau_{\rm relax}\) (\(\Delta t_{\rm div}< 0.2\)), it reduces the wall time needed by more than 50\%. For cell proliferation occurring on time scales longer than the relaxation time, the reduction increases further, with wall times for \(\Delta t_{\rm div}= 1.0\) being reduced by nearly 70\% and wall times for \(\Delta t_{\rm div}= 5.0\) by almost 90\%, decreasing the absolute run time from 75 minutes with the fixed time stepping to less than 9 minutes. The comparison in \eqref{eq:workdet} and \eqref{eq:workstoch} shows that similar results can be expected when the proliferations occur randomly.

Algorithms II (SRFES) and III (MRFE), which have a significantly higher cost per time step than Algorithm I due to calculation of the stability bound, are beneficial compared to fixed time stepping for \(\Delta t_{\rm div} > 0.5\). As expected, the more frequent cell proliferation, the more advantageous the use of the multirate Algorithm III is compared to the single rate Algorithm II. For highly frequent cell divisions (\(\Delta t_{\rm div} < 0.5\)), however, fixed time stepping is more efficient than both algorithms. The exact limit of this trade-off is implementation dependent and could potentially be shifted in favor of Algorithm III by further optimization of the code. For \(\Delta t_{\rm div} > 2.0\), Algorithm II gains an advantage over Algorithm III as the system spends long periods of time limited by stability where the algorithms choose the same time step \(\Delta t_s\), but Algorithm III requires additional overhead. This means there is an interval of \(\Delta t_{\rm div}\) values for which the multirate method, Algorithm III, is better than the single rate method Algorithm II. The gain of using Algorithm II over fixed time stepping increases with larger values of \(\Delta t_{\rm div}\), from a reduction of 20\% for \(\Delta t_{\rm div} = 1.0\) to a reduction of 70\% for \(\Delta t_{\rm div} = 5.0\).  

Using Algorithm IV (SRBE) is inefficient for high cell division frequencies (small values of \(\Delta t_{\rm div}\)). When cell divisions are sufficiently rare and the cell system is approaching steady state, however, it becomes beneficial compared to fixed time stepping. For our implementation and experimental setup this is the case when \(\Delta t_{\rm div} = 3.0\). 
For \(\Delta t_{\rm div} = 5.0\), i.e.\ there is one cell division every five multiples of the relaxation time \(\tau_{\rm relax}\), simulation with Algorithm IV is faster than with the multirate Algorithm III (MRFE). Extrapolating the data for \(\Delta t_{\rm div}\), one can expect Algorithm IV to be competitive compared to Algorithm II and possibly even Algorithm I for \(\Delta t_{\rm div}\) values larger than 5.0. 

\begin{table}[tbp]
\small
\renewcommand{\arraystretch}{1.25}
\centering
\begin{tabular}{c|c|r|r|r|r|r}
& & \multicolumn{5}{c}{Wall time [s]} \\
\(\Delta t_{div}\) & \(T\) &  Alg.\ I (SRFE) & Alg.\ II (SRFES) & Alg.\ III (MRFE) & Alg.\ IV (SRBE) & \(\Delta t_{fixed}\)   \\
\hline
0.1 & 1.0 & 83.58& 216.34&  106.77& 551.63& 86.69\\
0.5 & 5.0 & 139.07& 442.93&  332.43& 966.55&422.95\\
1.0 & 10.0 &  167.11& 704.39&  601.08& 1273.63& 892.16\\
5.0 & 50.0 &  534.82& 1340.91& 2826.74&  2016.27&4498.78\\
\end{tabular}
\caption{Simulation end times \(T\) and absolute wall times for Algorithms I-IV and the fixed time stepping runs for the different \(\Delta t_{\rm div}\) values used in the linear tissue growth experiment. The initial spheroid size is chosen as \(N_0 = 2197\). Simulation end times \(T\) are chosen to let 10 cell divisions at a frequency of \(1/\Delta t_{\rm div}\) occur. The fixed time step size \(\Delta t_{\rm fixed}= 0.0078\) is given by the minimal possible time step size after proliferation in \eqref{eq:Aprolifrestr}. Absolute wall times are averaged over four repetitions and are given in seconds rounded to two decimal values. }
\label{tab:linear_growth}
\end{table}

We can conclude that adaptive time stepping using Algorithm I in Figure \ref{fig:linear_tissue_growth} is universally favorable. This is also the conclusion from Figure \ref{fig:benchmark_single_prolif_spheroid}(c). If accuracy is of higher concern and Algorithm I is insufficient in that regard, see Figure \ref{fig:convergence_error_spheroid}, using Algorithm II still results in significant gains over standard fixed time stepping unless cell proliferation is very frequent. This may be the case in large cell aggregations and the fixed time step method would be the preferred choice. Then the constant time step should be determined after the first proliferation by \eqref{eq:forwardEdt}.

\section{Summary and Conclusions}\label{sec:concl}
Center-based models for the forces between biological cells are analyzed and numerical methods for simulation of the models are proposed and compared. The analytical equations form a system of ODEs of first order. 
The forces are obtained as the gradient of a potential. The numerical methods are developed for center-based models for cell simulations but are generally applicable to gradient systems of ODEs \eqref{eq:poteq}, for example vertex-based models of the cell forces, and to other force functions than \eqref{eq:cubic_force}.

The equations are approximated by the forward Euler and backward Euler methods of first order accuracy. 
The backward Euler method has similar stability properties as the analytical problem with a solution converging to an equilibrium which is constant in time.
The equilibrium is not unique but depends on the initial conditions.
If the solution reaches an equilibrium solution it remains there in the analytical solution and the Euler solutions.
The center of gravity is constant in the analytical and numerical  solutions. 
The eigenvalues of the Jacobian matrix of the forces are real. Close to an equilibrium point they are non-positive, but at least after cell proliferation there are positive eigenvalues.

The local error in a time step is estimated. An error bound or a stability bound  determines the time step for the forward Euler method. 
The eigenvalues of the Jacobian matrix are estimated for the stability bound. Only the local error bounds the time step for the backward Euler method. 
A multirate time stepping method is proposed. The positions of some cells are integrated by shorter time steps while longer time steps are taken for the majority of the cells. A bound on the global discretization error is derived for the forward Euler method using the bound on the local error.
The performance of four methods with time step adaptation are compared in numerical experiments: the single rate forward Euler method with and without a stability check (SRFES and SRFE), the multirate forward Euler method with local time steps (MRFE), and the single rate backward Euler method (SRBE).  

We have shown that Algorithm I, the single rate forward Euler method, leads to a strong decrease in computational cost --- a reduction of 70\% to 90\% in terms of wall time for the simulation of a linearly growing spheroid--- compared to a fixed time stepping forward Euler method that uses the same initial time step (restricted by strong repulsive forces between daughter cells after proliferation) at the cost of a slightly decreased accuracy. 
There are oscillations in the time step sequence and the global errors but at a low level. We hence propose the SRFE algorithm for general use since it is easy to implement, eliminates the need to manually  determine the correct time step size, and offers a more intuitive absolute error in each time step on the spatial scale instead. 

If accuracy is of a higher concern, our results confirm that unless cell division on the population level happens very much faster than the mechanical relaxation between two individual daughter cells, both Algorithms II (SRFES) and III (MRFE) 
beat the fixed time stepping algorithm in terms of wall time and computational cost. Note that our findings are mostly independent of the exact system size.
The exact value of how much time needs to pass between cell division events at the population level so that the overhead of either algorithm pays off is of course implementation dependent. For our results it lies at roughly \(0.5\tau_\text{relax}\). Typical cell cycle durations depend strongly on the cell type, for human cells 22-24 hours are common \cite{cooper2007cell}. On the other hand the mechanical relaxation time can be considered to be on the order of a few minutes or even less. Consequently, there will be a range of population sizes for which the cell division frequencies we consider here are applicable.  
An alternative to the Euler methods is to upgrade the order of accuracy with an explicit Runge-Kutta method of order two. Theorem 2 is still valid and the error estimate in \eqref{eq:forwardEerrest} should be changed to an estimate of $\d^3\fatx/\d t^3$.


From Figures \ref{fig:2cells_dt}, \ref{fig:prolif_within_spheroid_dt} and \ref{fig:dependence_n_cells} we see that the SRBE method in Algorithm IV is able to take longer steps while the other three algorithms are restricted by stability concerns, which for both configurations (i) and (ii) occurs shortly after the relaxation time. Nevertheless, the SRBE method achieves its improved stability properties at an increased computational cost per step. 
It is therefore not straightforward to know if there is a gain using the SRBE method in terms of a reduced wall clock time for the simulation as discussed in Section \ref{sec:SRBE}. In fact, in \cite{Atwell2016,mathias2022cbmos} it was shown for different population configurations that the performance of a fixed time stepping backward Euler method was not better than a fixed time stepping forward Euler method. 
In this study, we showed that using adaptive time stepping for the backward Euler method makes it more efficient than a fixed time stepping forward Euler method when cell proliferation is sufficiently rare and the system spends long duration limited by stability. However, SRBE is less efficient than the adaptive time stepping methods based on the forward Euler method for the cell proliferation frequencies we consider. 

The efficiency of the algorithms depends to some extent on the implementation of them. There is a large difference in computational cost between the algorithms that rely on explicit calculation of the stability bound (SRFES and MRFE) and SRFE that does not. It is expensive to assemble the Jacobian necessary to calculate or estimate its eigenvalues which determine the stability bound. In our Python implementation we assemble the complete Jacobian in order to benefit as much as possible from NumPy's routines. It is, however, possible ---and highly advisable for a compiled language such as C/C++--- to apply Gershgorin's estimate by only calculating relative parts of the Jacobian on the fly. Similarly, the computational cost can be reduced by not recalculating the Jacobian (or Gershgorin's estimate) unnecessarily once the system is in a steady state, but ``freezing'' it until the next cell proliferation happens. These modifications would increase the gain of SRFES and MRFE compared to fixed time stepping. 

Other possibilities to extend our work include its application to a variant of the center-based model where the neighborhood definition is based on the Voronoi tesselation given by the cell midpoints \cite{VLPJD15, Meineke2001}. 
In this case, cells only interact with direct neighbors sharing an edge in the Voronoi tesselation. Movement of the cells and resulting rearrangement of neighbors potentially leads to discontinuous changes in the right hand side of the ODE system but this should not cause any problems for the adaptive algorithms based on the Euler method. 
\section*{Acknowledgment}\label{sec:ack}
The authors would like to thank Andreas Hellander and Adrien Coulier for fruitful discussions around the content of this article and the implementation of the algorithms and for comments on the manuscript. 
This work has received funding from the NIH under
grant no. NIH/2R01EB014877-04A1 and from the eSSENCE strategic initiatives on eScience. The funders had no role in the design of the study, data collection, data analysis, interpretation of results, or writing of the manuscript.
Numerical experiments were performed on the Rackham compute resources provided through the Uppsala Multidisciplinary Centre for Advanced Computational Science (UPPMAX) within the Project SNIC 2021-22-607.

\section*{Declarations of interest}\label{sec:dec}
None.

\bibliographystyle{plain}
\bibliography{main}

\end{document}